\documentstyle[epsf,amssymb,12pt]{article}
\newcommand{\aalpha}{{\overline{\alpha}}}
\title{Dual Teichm\" uller spaces.}
\author{V.V.Fock \\ {\em Institute of
Theoretical and Experimental Physics,}\\
{\em 117259, B.Cheremushkinskaya 25, Moscow, Russia.}\\ { \em E-mail
{\tt fock@mpim-bonn.mpg.de}} }
\date{25 December 1997}

\begin{document}
\maketitle

\begin{abstract}
We describe two spaces related to Riemann surfaces --- the Teichm\" uller
space of decorated surfaces and the Teichm\" uller space of surfaces with holes.
We introduce simple explicit coordinates on them. Using these coordinates we
demonstrate the relation of these spaces to the spaces of measured
laminations, compute Weil-Petersson forms, mapping class group action and
study properties of lamination length function. Finally we use the developed
technique to construct a noncommutative deformation of the space of functions
on the Teichm\" uller spaces and define a class of unitary projective mapping
class group representations (conjecturally a modular functor). One can
interpret the latter construction as quantisation of 3D or 2D Liouville
gravity. Some theorems concerning Markov numbers as well as Virasoro orbits
are given as a by-product.
\end{abstract}

\section{Introduction.}

The main philosophical aim of the paper is to {\em formulate} two problems
concerning Teichm\" uller spaces of Riemann surfaces with holes.

The first problem is to describe explicitly a kind of Fourier transform
between the spaces of functions on two slightly different versions of
Teichm\" uller spaces.

The second problem is to deform (quantise) the algebra of
functions on Teichm\" uller space in a direction prescribed by the
Weil-Petersson Poisson bracket and compatible w.r.t. homotopy classes of
mappings between surfaces. In particular w.r.t. the mapping
class group action.

The solution for the second problem is given in this version of the text,
however we do not give here detailed proofs and examples which is postponed
to a subsequent paper. Concerning the first problem, we just try to
give arguments in favor of existence of a solution and emphasise its
importance.

The technique used in the article unifies the Thurston approach to Riemann
surfaces, such as measured laminations on the one hand, and mathematical physics such as modular functors of conformal field theory on the other. One
metaresult important for us was the construction of a bridge between these
domains. 

  However the article does not contain philosophical discussions except for
a few remarks.  In the main part of the text we give definitions and prove
theorems (which can be considered as preparatory in the spirit of the
problems described above, but we hope that they have some independent
interest as well). The main ones of them are:

  We define explicitly simple global coordinates on the Teichm\" uller spaces
of Riemann surfaces with holes. We describe Penner's coordinates on the
decorated Teichm\" uller spaces. We give explicit formulae for the action of
the mapping class group as well as for the Weil-Petersson Poisson bracket on
the former space and for the Weil-Petersson degenerate symplectic structure
for the latter one. We show using the coordinates that these spaces have 
natural "scaling" limits to the space of measured laminations with closed
(resp.  compact) support. We give an elementary proof of continuity of the
lamination length function and the lamination intersection number. As a
by-product of the latter statement we prove some continuity theorem
concerning Markov numbers.  We show also explicitly compatibility of the
length and the intersection index functions at the limit when Teichm\"uller
spaces go to the respective laminations. We describe Bers's coordinates on the
simplest Virasoro orbit and as another by-product we compute the Kirillov--Kostant Poisson bracket in terms of these coordinates.
 Finally we give an explicit construction of noncommutative
deformation of the space of functions on the Teichm\' uller spaces of Riemann
surfaces with holes depending on a quantisation parameter $\hbar$ and show an
amazing symmetry between deformations corresponding to the parameters $\hbar$
and $\frac{1}{\hbar}$.

We tried to make the paper to be self-contained and available for a wide class
of readers. Therefore we have included many known results. Some of them are
provided with a few line proofs in the spirit of the paper. The other
proofs are left for the interested reader as easy exercises. Some slightly
more complicated results in the two final sections are provided with
references to a proof.

  For a nonrigorously minded reader we remark that the quantisation procedure gives Hilbert spaces which can be interpreted either as the space of conformal
blocks  of Liouville gravity theory in two Euclidean dimensions or as the
space of states in 3D quantum gravity, as it follows form the ideas of \cite{Witten}, \cite{VV}. The $\hbar \leftrightarrow 1/\hbar$ symmetry is rather similar to the one observed in \cite{ZZ}. However we omit (except for a few
sentences at the end of some sections) the discussion of this point of view
here since we can hardly imagine arguments making something more out of this statement than just a definition.

\section{Graphs and surfaces.}
  In this section we shall give a brief description of relations between
surfaces and fat graphs. These relations exist only for surfaces with the
number of holes $s \ge 1$, genus $g \ge 1$ or with at least 3 holes and genus
$0$. Such surfaces will be called {\em hyperbolic}. All graphs considered in thesequel are supposed to be finite.

 Recall that a {\em fat graph} is an unoriented graph s.t. for each vertex the cyclic order of ends of edges incident to the vertex is given.

One can imagine a fat graph as a graph with edges being narrow bands. (It
is where the attribute {\em fat} comes from.) A graph drawn on an oriented surface
acquires a fat graph structure given by, say, a counterclockwise ordering of
the ends of edges at each vertex.

Let us say that an oriented path on a fat graph {\em turns left} at a vertex if we come to the vertex along the end of an edge which is precedent to the one we come out w.r.t. the cyclic order. 

 On a fat graph one can well define a distinguished set of closed paths called
{\em faces}. A face is a path s.t. being oriented it turns always left at each verex or always right. 

Denote by $V(\Gamma),E(\Gamma)$ and $F(\Gamma)$ the sets of vertices, edges
and faces of $\Gamma$, respectively.

 We can obtain a smoothable surface $S_0(\Gamma)$ from a fat graph $\Gamma$ by
taking a disk for each face and gluing its boundary to the graph along the
face. If we have taken the set of edges and vertices of a polyhedron in ${\Bbb
R}^3$ as a graph with a natural fat structure, we recover by this procedure the
original polyhedron, the faces of the graph being correspondent to the faces of
the polyhedron. And this is the reason to use the term {\em face} in this
context. We can use annuli instead of disks and get a surface $S(\Gamma)$ with
boundary. The surface $S(\Gamma)$ can be obviously retracted onto $\Gamma$.

 One can give a purely combinatorial description of a fat graph. Let
$EE(\Gamma)$ be the set of ends of edges of a fat graph $\Gamma$ (or what is the
same, the set of oriented edges).
Define an involution $s_1$ acting on $EE(\Gamma)$ which maps an end of an
edge to the opposite end of the same edge (resp. reverse the orientation of
the edge).  The fat structure induces a permutation $s_0$ of the same set
which maps an end of an edge to the next end of an edge w.r.t. the cyclic order
at the corresponding vertex. It is obvious that this gives a one-to-one
correspondence between fat graphs and pairs of permutations
$(s_0,s_1)$ on finite sets, s.t. $s_1$ is an involution without fixed points.
In these terms faces correspond to orbits of $s_0^{-1}s_1$ in $EE(\Gamma)$.

  Introduce some notation useful for the sequel. Let $\aalpha \in EE(\Gamma)$
be an end of an edge. Then let
$$
\aalpha(1) = s_0 \aalpha;~~
\aalpha(2) = s_0^{-1} s_1 \aalpha;~~
$$ $$
\aalpha(3) = s_0 s_1 \aalpha;~~
\aalpha(4) = s_0^{-1} \aalpha.
$$
\setlength{\unitlength}{0.025em}
\begin{picture}(0,0)(-1100,600)
\thicklines
\put( 40,640){\vector( 1, 2){ 40}}
\put(260,640){\vector(-1, 2){ 40}}
\put(260,800){\vector(-1,-2){ 40}}
\put( 40,800){\vector( 1,-2){ 40}}
\put( 80,720){\vector( 1, 0){140}}
\put(170,730){\makebox(0,0)[rb]{$\aalpha$}}
\put( 60,640){\makebox(0,0)[lb]{$\aalpha(2)$}}
\put( 55,780){\makebox(0,0)[lb]{$\aalpha(3)$}}
\put(245,780){\makebox(0,0)[rb]{$\aalpha(4)$}}
\put(245,640){\makebox(0,0)[rb]{$\aalpha(1)$}}
\end{picture}
We say that $\aalpha \in \gamma$, where $\aalpha \in EE(\Gamma)$ and $\gamma \in
F(\Gamma)$ if the orientation of $\aalpha$ agrees with the counterclockwise
orientation of $\gamma$. Denote by $\alpha$ the edge from $E(\Gamma)$ corresponding to $\aalpha \in EE(\Gamma)$.

 Note that for a fat graph $\Gamma$ one can define the dual graph
$\Gamma^\vee$ with vertices, edges and faces replaced by faces, edges and
vertices of $\Gamma$, respectively.  In terms of permutations the dual graph
corresponds to the pair $(s_2,s_1)$, where $s_2 = s_0^{-1} s_1,$ acting on
the same set.

 Now show that any hyperbolic surface $S$ can be obtained as $S(\Gamma)$ for
some fat graph ${\Gamma}$.  The construction becomes more transparent if we
imagine the holes having zero size, i.e., just as punctures.  Any surface can
be cut into topologically trivial pieces by a number of curves going from
puncture to puncture which are self- and mutually nonintersecting, mutually
homotopically nonequivalent and nonshrinkable to punctures. We can always
take enough curves to make the pieces simply connected.  The resulting set of
curves is a graph ${\Gamma}^\vee$ with vertices at the punctures and a
natural fat structure given by the orientation of the surface.  The desired
graph ${\Gamma}$ is dual to this one.  One can easily check by drawing
pictures that the surface $S(\Gamma)$ is isomorphic to the surface $S$ we
have started with.  If we now take one point inside each piece and for each
cut draw a segment with ends in the chosen points intersecting only this cut
we obtain the graph $\Gamma$ together with the homotopy class of its embedding
into $S$.

  A maximal system of such curves (it always exists) cuts the surface into
triangles and the corresponding graphs turn out to have three ends of edges
incident to each vertex (3-valent graphs). In the sequel we mostly consider
such kind of graphs. One can easily check by computing the Euler
characteristics that a 3-valent graph $\Gamma$ has $6g-6+3s$ edges $4g-4+2s$
vertices and $s$ faces, where $g$ and $s$ are the genus and the number of
holes of the surface $S(\Gamma)$, respectively.

Denote by $\Gamma(S)$ the set of graphs, corresponding to a given surface $S$
and by $\Gamma_0(S)$ the set of graphs together with their embeddings into
$S$ considered up to homeomorphisms of $S$ homotopy equivalent to the
identity. The former set is finite and the latter is obviously infinite.
The mapping class group ${\cal D}(S)$ acts naturally on $\Gamma_0(S)$ with
$\Gamma(S)$ as a quotient.  The subset of $\Gamma_0(S)$ (resp.
$\Gamma(S))$ consisting of three-valent graphs is denoted by $\Gamma_0^3(S)$
(resp. $\Gamma^3(S)).$ It is obviously stable w.r.t. ${\cal D}(S)$.

\noindent
\begin{minipage}{9.5cm}
{ ~~~~There exists a natural operation called {\em flip} or {\em
Whitehead move} which makes one element of  $\Gamma_0^3(S)$ from another.
Consider an edge of the dual graph which bounds two triangles forming
together a quadrilateral.  (On the original graph this condition means that
the ends of the edge do not coincide). Remove this edge and replace it by the
second diagonal of this quadrilateral. On the original graph this operation
means shrinking an edge and then blowing it up in another direction as
shown on the picture \ref{jflip}. }
\end{minipage}
\hfill
\begin{minipage}{6cm} 
\unitlength 0.33mm
\linethickness{0.4pt}
\vspace{0.5cm}
\centerline{
\begin{picture}(120,100)(15,15)
\multiput(80.00,120.00)(0.12,0.12){167}{\line(0,1){0.12}}
\multiput(100.00,140.00)(0.12,-0.12){167}{\line(0,-1){0.12}}
\multiput(120.00,120.00)(-0.12,-0.12){167}{\line(0,-1){0.12}}
\multiput(100.00,100.00)(-0.12,0.12){167}{\line(0,1){0.12}}
\put(80.00,120.00){\line(1,0){40.00}}
\multiput(30.00,140.00)(0.12,-0.12){167}{\line(0,-1){0.12}}
\multiput(50.00,120.00)(-0.12,-0.12){167}{\line(0,-1){0.12}}
\multiput(30.00,100.00)(-0.12,0.12){167}{\line(0,1){0.12}}
\multiput(10.00,120.00)(0.12,0.12){167}{\line(0,1){0.12}}
\put(30.00,140.00){\line(0,-1){40.00}}
\put(75.00,120.00){\vector(1,0){0.2}}
\put(55.00,120.00){\line(1,0){20.00}}
\put(-10.00,120.00){\makebox(0,0)[lc]{$\Gamma^\vee$}}
\multiput(10.00,80.00)(0.12,-0.24){84}{\line(0,-1){0.24}}
\put(20.00,60.00){\line(1,0){20.00}}
\multiput(40.00,60.00)(0.12,0.24){84}{\line(0,1){0.24}}
\multiput(50.00,40.00)(-0.12,0.24){84}{\line(0,1){0.24}}
\multiput(10.00,40.00)(0.12,0.24){84}{\line(0,1){0.24}}
\multiput(80.00,40.00)(0.24,0.12){84}{\line(1,0){0.24}}
\put(100.00,50.00){\line(0,1){20.00}}
\multiput(100.00,70.00)(-0.24,0.12){84}{\line(-1,0){0.24}}
\multiput(120.00,80.00)(-0.24,-0.12){84}{\line(-1,0){0.24}}
\multiput(120.00,40.00)(-0.24,0.12){84}{\line(-1,0){0.24}}
\put(75.00,60.00){\vector(1,0){0.2}}
\put(55.00,60.00){\line(1,0){20.00}}
\put(-10.00,60.00){\makebox(0,0)[lc]{$\Gamma$}}
\end{picture}}
\refstepcounter{equation}\label{jflip}
\centerline{Fig. \arabic{equation}}
\end{minipage}

If a graph $\Gamma \in \Gamma^3(S)$ has a symmetry, it acts obviously on the
corresponding elements of $\Gamma_0^3(S)$.
One can show that any two elements of $\Gamma_0^3(S)$ are connected by
a sequence of flips and graph symmetries. In particular, the action of any
element of the mapping class group ${\cal D}(S)$ on $\Gamma_0^3(S)$ can be
represented by a sequence of flips and symmetries. In more scientific words
$\Gamma_0^3(S)$ and sequences of flips and symmetries constitute a groupoid
containing the mapping class group as the greatest subgroup.

Note that if $\sigma: \tilde{S} \rightarrow S$ is an unramified $N$-fold
covering then one can obviously construct a graph corresponding to $\tilde{S}$
starting from a graph $\Gamma$ corresponding to $S$. (This graph
$\sigma^*\Gamma$ is just the full inverse image of $\Gamma$ in
$\tilde{S}$.  $\sigma^*\Gamma$ has $N \sharp E(\Gamma)$ edges, $N \sharp
V(\Gamma)$ vertices and $\sum_{\gamma \in F(\Gamma}) O(\gamma)$ faces.  Here
$O(\gamma)$ is the number of orbits of the covering monodromy around the face
$\gamma$.  There is a natural mapping from the edges, vertices and faces of
$\sigma^*\Gamma$ to the edges, vertices and faces of $\Gamma$, respectively,
which we shall denote by the same letter $\sigma$.

  The mapping class group ${\cal D}(S)$ obviously acts on the set of
unramified $N$-fold coverings of $S$. A stabiliser of a covering $\sigma$ in
${\cal D}(S)$ we call a {\em congruence subgroup w.r.t. $\sigma$} and denote
by ${\cal D}(S,\sigma)$. ${\cal D}(S,\sigma)$ is obviously a finite index
subgroup in ${\cal D}(S)$

 Call a three-valent fat graph {\em regular} if it
has no edges with coinciding ends, any two edges have no more than one
common vertex and any edge separates two different faces.  Not all surfaces
can be represented by regular graphs.  The only reason to introduce this
class of graphs is because usually all constructions and formulae are more
simple for them.  However any nonregular graph can be covered by a regular
one, and usually one can easily derive formulae for nonregular graphs
starting from those for regular graphs by passing to such a covering.

\section{Laminations.}

 Taking into account that the reader may not be familiar with the Thurston's
notion of a measured lamination \cite{Thurston}, we are going to give all
definitions here in the form, which is almost equivalent to the original one
(the only difference is in the treatment of the holes and punctures), but more
convenient for us. The construction of coordinates on the space of laminations
we are going to describe is a slight modification of Thurston's "train tracks"
(\cite{Thurston}, section 9 and \cite{Trtr}).

It seems worth mentioning here, that the definitions of measured
laminations are very similar to the definitions of the singular homology
groups, and is in a sense an unoriented version of the latter ones.

There are two different ways to define the notion of measured
laminations for surfaces with boundary, which are analogous to the definition
of homology group with compact and closed support, respectively.

\subsection{Bounded measured laminations.}

\paragraph{Definition.}

{\em Rational bounded measured lamination} on a 2-dimensional surface is a
homotopy class of a collection of finite number of self- and mutually
nonintersecting unoriented closed curves  with rational weights and subject
to the following conditions and  equivalence relation.

1. Weights of all curves are positive, unless a curve surrounds a hole.

2. A lamination containing a curve of weight zero is considered to be
equivalent to the lamination with this curve removed.

3. A lamination containing two homotopy equivalent curves with weights $a$
and $b$ is considered as equivalent to the lamination with one of these curves
removed and with the weight $a+b$ on the other.

 The set of all rational bounded laminations on a given surface $S$ is denoted
by ${\Bbb Q}{\cal L}^d(S).$ This space has a natural subset, consisting of
laminations with integral weights.  The set of such laminations is denoted by
${\Bbb Z}{\cal L}^d(S)$. Denote by ${\Bbb Q}{\cal L}(S) \in {\Bbb Q}{\cal
L}^d(S)$ and ${\Bbb Z}{\cal L}(S) \in {\Bbb Z}{\cal L}^d(S)$ the subspaces,
consisting of laminations without curves surrounding holes.

 {\em Remark.} Any rational bounded measured lamination can be represented by a
collection of $3g-3+s$ curves. Any integral lamination can be represented by a
finite collection of curves with weights $+1$ or $-1$ on some curves surrounding
holes.

\paragraph{Construction of coordinates.}

Suppose we are given a three-valent fat graph $\Gamma \in \Gamma_0^3(S)$.  We
are going to assign, for a given lamination, rational numbers on edges
and show, that these numbers are good coordinates on the space of laminations.

\noindent\begin{minipage}[c]{9cm} 
~~~Retract the lamination to the graph in such a way, that each curve retracts to a path without folds on edges of the graph, and no curve goes along
an edge and then back, without visiting another edge.  Assign to each edge
$\alpha$ the sum of weights of curves, going through it (fig.  \ref{blam}).
The collection of these numbers, one for each edge of $\Gamma$, is the desired
set of coordinates.
\end{minipage}
\begin{minipage}{6cm}
\epsfxsize4.5cm
\centerline{\epsfbox{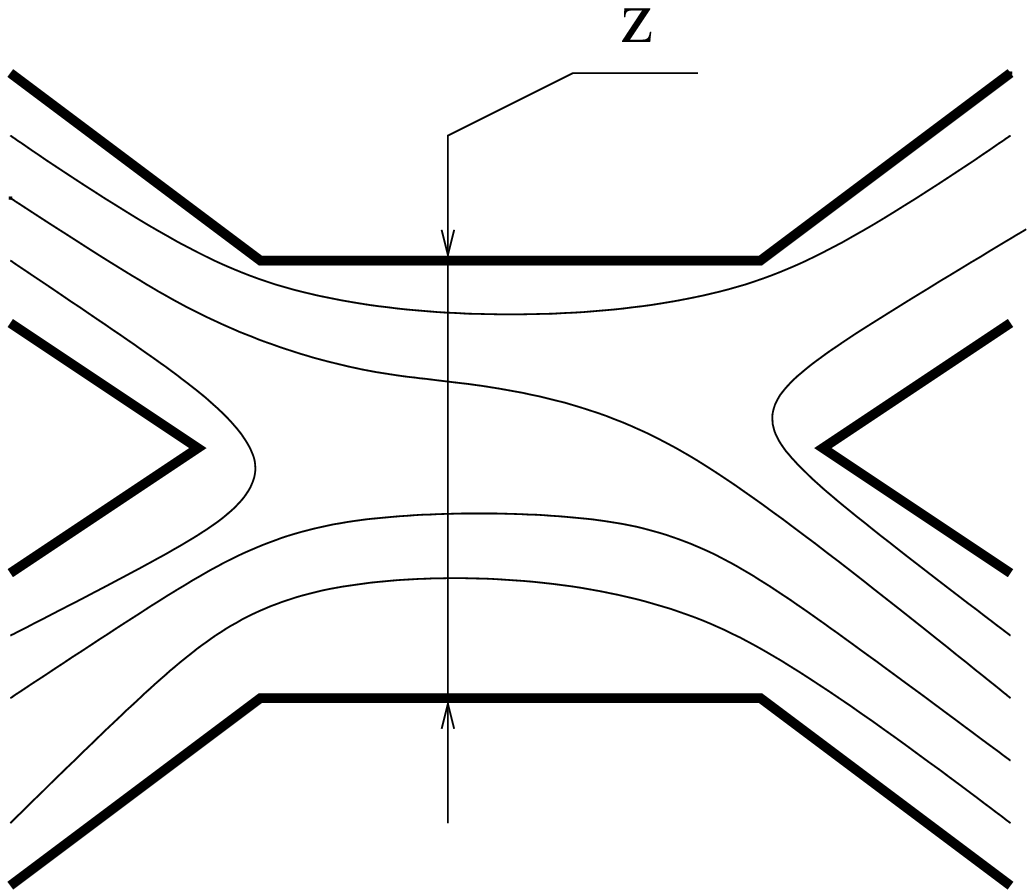}}
\refstepcounter{equation}\label{blam}
\centerline{Fig.\arabic{equation}}
\end{minipage}

\paragraph{Reconstruction.}
Now we need to prove that these numbers  are coordinates indeed. For this
purpose we just describe an inverse construction which gives a lamination,
starting from numbers on edges.

 First of all note, that if we are able to reconstruct a
lamination, corresponding to a set of numbers $\{z_\alpha\}$ , we can do it as
well for the set $\{az_\alpha\}$ and $\{z_\alpha+b\}$ for any rational $a \ge
0$ and $b$. Indeed, multiplication of all numbers by $a$ can be achieved by
multiplication of all weights by $a$ and adding $b$ is obtained by adding
loops with weight $b/2$ around each hole. Therefore, we can use these
possibilities  to reduce our problem to the case when $\{z_\alpha\}$ are
positive integers and any three numbers on edges incident to each vertex $z_1,
z_2, z_3$ satisfy triangle and parity conditions
\begin{eqnarray}
\vert z_1-z_2 \vert \le z_3 \le  z_1+z_2  \\  z_1+z_2+z_3  \hbox{ is even}
\end{eqnarray}

 Now the reconstruction of the lamination is almost obvious. Draw
$z_\alpha$ lines on the $\alpha-$th edge and connect these lines at vertices in
a nonintersecting way (fig.\ref{vertex}), what can be done unambiguously.

\centerline{{\epsfxsize8\baselineskip\epsfbox{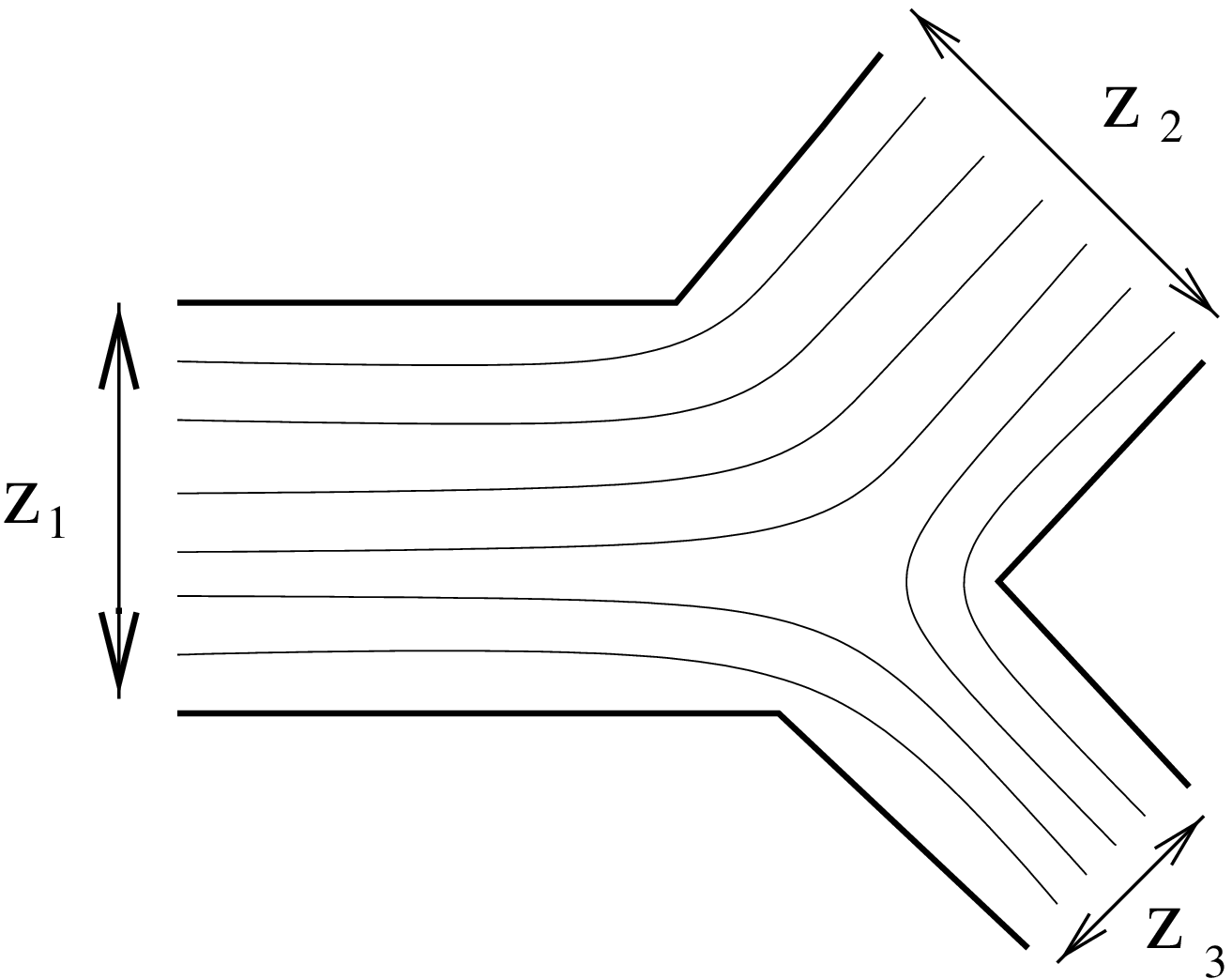}}}
\refstepcounter{equation}\label{vertex}
\centerline{Fig. \arabic{equation}} 

\paragraph{Graph change.}

 The constructed coordinates on the space of laminations is related to a
particular choice of the three-valent graph.  The following formulae describe
the change of coordinates under a flip of an edge of the graph.

\vspace{5mm}
\centerline{
\setlength{\unitlength}{0.17mm}%
\begin{picture}(500,275)(20,485)
\thicklines
\put(280,700){\line( 1,-2){ 40}}
\put(320,620){\line( 1, 0){280}}
\put(600,620){\line( 1, 2){ 40}}
\put(600,620){\line( 1,-2){ 40}}
\put(320,620){\line(-1,-2){ 40}}
\thinlines
\put(180,620){\vector(-1, 0){  0}}
\put(180,620){\vector( 1, 0){ 55}}
\thicklines
\put(100,540){\line( 2,-1){ 80}}
\put(100,700){\line( 0,-1){160}}
\put(100,540){\line(-2,-1){ 80}}
\put( 20,740){\line( 2,-1){ 80}}
\put(100,700){\line( 2, 1){ 80}}
\put( 40,745){\makebox(0,0)[lb]{$A$}}
\put(160,745){\makebox(0,0)[rb]{$B$}}
\put(120,620){\makebox(0,0)[lc]{$Z$}}
\put(160,485){\makebox(0,0)[rt]{$C$}}
\put( 40,485){\makebox(0,0)[lt]{$D$}}
\put(300,540){\makebox(0,0)[lt]{$D$}}
\put(620,540){\makebox(0,0)[rt]{$C$}}
\put(620,690){\makebox(0,0)[rb]{$B$}}
\put(300,690){\makebox(0,0)[lb]{$A$}}
\put(470,640){\makebox(0,0)[cb]{$\max(A+C,B+D)-Z$}}
\end{picture}}
\refstepcounter{equation}
\centerline{Fig. \arabic{equation}}
\label{dflipl}
\vspace{5mm}
(Only the changing part of the graph is shown here, the numbers on the other edges
remain unchanged.)

\subsection{Unbounded measured laminations}
\paragraph{Definition}

{\em Rational unbounded measured lamination} on a 2-dimensional surface with
boundary is a homotopy class of a collection of finite number of
nonselfintersecting and pairwise nonintersecting curves either closed or
connecting two boundary components (possibly coinciding) with positive
rational weights assigned to each curve and subject to the following
equivalence relations:

 1. A lamination, containing a curve retractable to a boundary component  is
equivalent to the lamination with this curve removed.

2. A lamination containing a curve of zero weight is considered to be
equivalent to the lamination with this curve removed.

 3. A lamination containing two homotopy equivalent curves of weights $a$ and
$b$, respectively, is equivalent to the lamination with one of these curves removed
and with the weight $a+b$ on the other.

 The set of all rational unbounded laminations on a given surface $S$ is
denoted by ${\Bbb Q}{\cal L}^h(S)$. This space has a natural subset,
representable by collections of curves with integral weights. This space is
denoted by ${\Bbb Z}{\cal L}^h(S)$.

 {\em Remark}.  Any rational unbounded measured lamination can be represented
by a collection of no more than $6g-6+2s$ curves (for Euler characteristics
reasons).  Any integral lamination can be represented by a finite collection of curves
with unit weights.

 For any given lamination, fix orientations of all boundary components but
those nonintersecting with curves of the lamination. Denote the
space of rational (resp. integral) laminations equipped with this additional
structure by ${\Bbb Q}{\cal L}^H$ (resp. ${\Bbb Z}{\cal L}^H$).

\paragraph{Construction of coordinates}

Suppose we are given a three-valent fat graph $\Gamma \in \Gamma_0^3(S)$.
We are going to assign for
a given element of the space ${\Bbb Q}{\cal L}^H$ a set of rational numbers on
edges, and show that these numbers are good global coordinates on this space.

 Straightforward retraction of an unbounded lamination onto $\Gamma$ is not
good because some curves may shrink to points or finite segments. To avoid this
problem, let us first rotate each oriented boundary component infinitely many
times in the direction prescribed by the orientation as shown on fig.
\ref{twist}.

\vspace{5mm}
\epsfxsize8cm
\centerline{\epsfbox{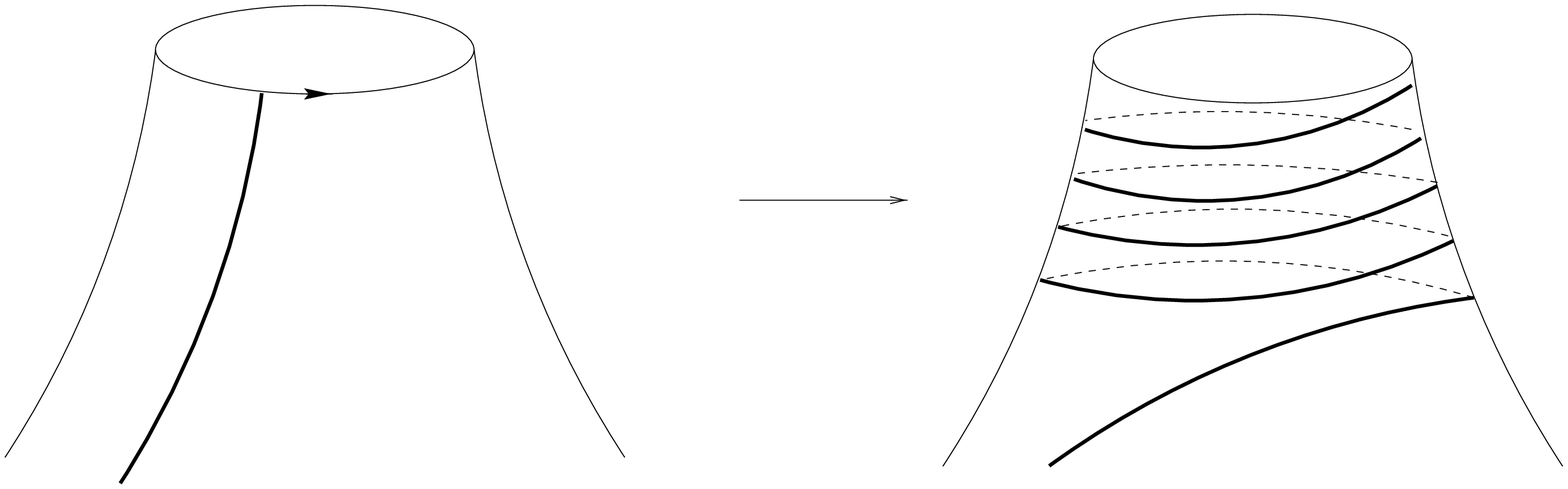}}
\refstepcounter{equation}\label{twist}

\vspace{5mm}
\centerline{Fig. \arabic{equation}}
\vspace{5mm}

\noindent\begin{minipage}{11cm}
~~~ The resulting lamination can be retracted on $\Gamma$ without
folds. Although we possibly get infinitely many curves, going through an edge.
Call that a curve {\em is right handed (resp., left handed)} in an edge if it turns left (resp., right) at both ends of the edge w.r.t. the motion along it from the center of the edge. Now assign to the edge the sum of weights of curves right handed in it with the sign plus or if this set is empty the sum of weights of left handed curves with the sign minus (fig. \ref{ublam}).
  
~~~The collection of these rational numbers, one for
each edge, is the desired coordinate system on ${\Bbb Q}{\cal L}^H$.
  
~~~Note that the number of right or left handed is always finite and
therefore the numbers assigned to the edges are correctly defined. Indeed,
consider the curves retracted on the graph. \end{minipage}
\hfill
\begin{minipage}{4cm}
\centerline{\epsfxsize3cm\epsfbox{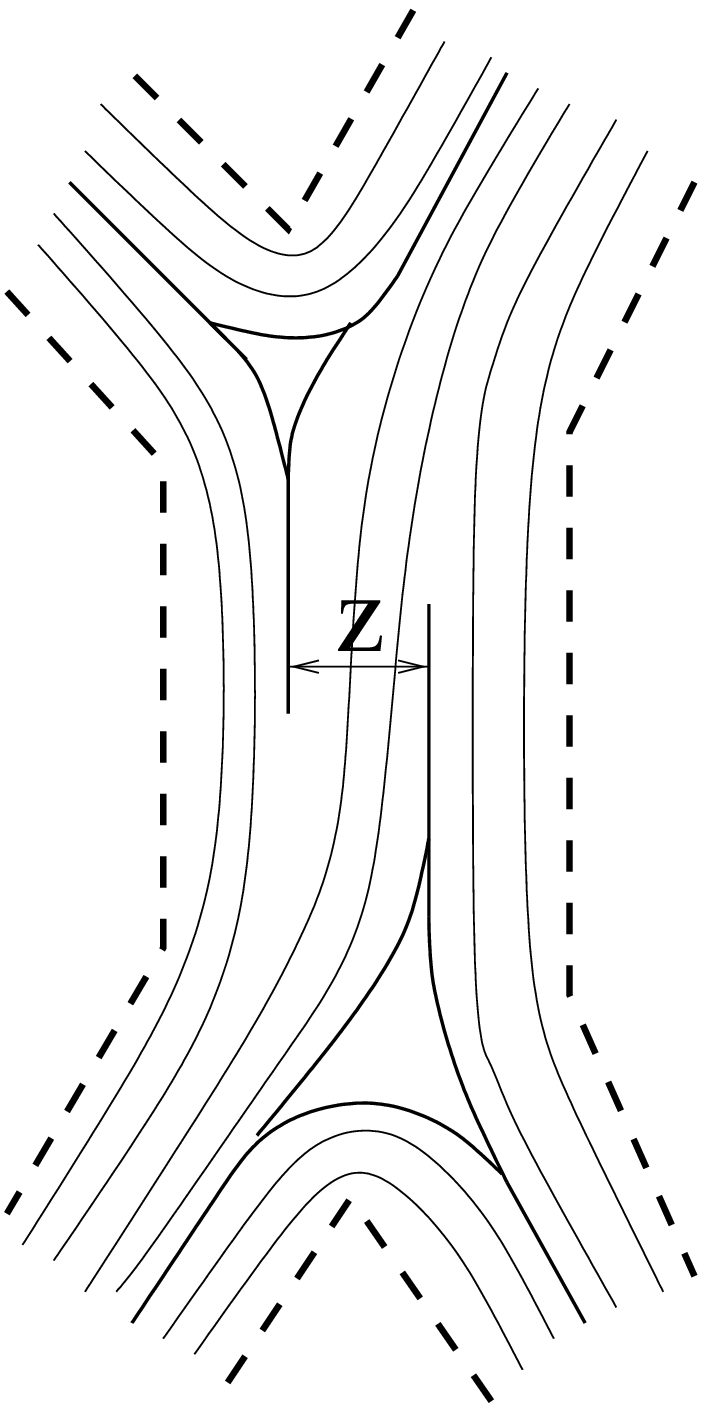}} 
\refstepcounter{equation}\label{ublam}
\centerline{Fig. \arabic{equation}}
\end{minipage}
\noindent We can mark a finite segment of
each nonclosed curve in such a way that each of two unmarked semiinfinite
rays goes only around a single face and therefore are never right or left handed. Therefore only the finite marked parts of curves contribute to
the numbers on the edges.

\paragraph{Reconstruction}

Now we need, as in the bounded case, to prove that these numbers are
indeed coordinates, what we shall do as well by describing an inverse
construction.  Note that if we are able to construct a lamination
corresponding to the set of numbers $\{u_\alpha\}$, we can equally do it for
the set $\{a u_\alpha\}$ for any rational $a \ge 0$.  Therefore we can reduce
our task to the case when all numbers on edges are integral.  Now draw
${\Bbb Z}$-infinitely many lines along each edge.  In order to connect these
lines at vertices we need to split them at each of the two ends into two
${\Bbb N}$-infinite bunches to connect them with the corresponding bunches of
the other edges.  Let us make it at the $\alpha$-th edge, such that $u_\alpha
\leq 0$ (resp.  $u_\alpha \ge 0$), in such a way, that the intersection of
the right (resp.  left) bunches at both ends of the edge consist of
$u_\alpha$ lines (resp.  $-u_\alpha$ lines).  Here the left and the right
side are considered from the centre of the edge toward the corresponding end.
The whole procedure is illustrated on fig.  \ref{ublam-}. The resulting
collection of curves may contain infinite number of curves surrounding holes,
which should be removed in accordance with the definition of an unbounded
lamination.

\vspace{5mm}
\centerline{{\epsfxsize16\baselineskip\epsfbox{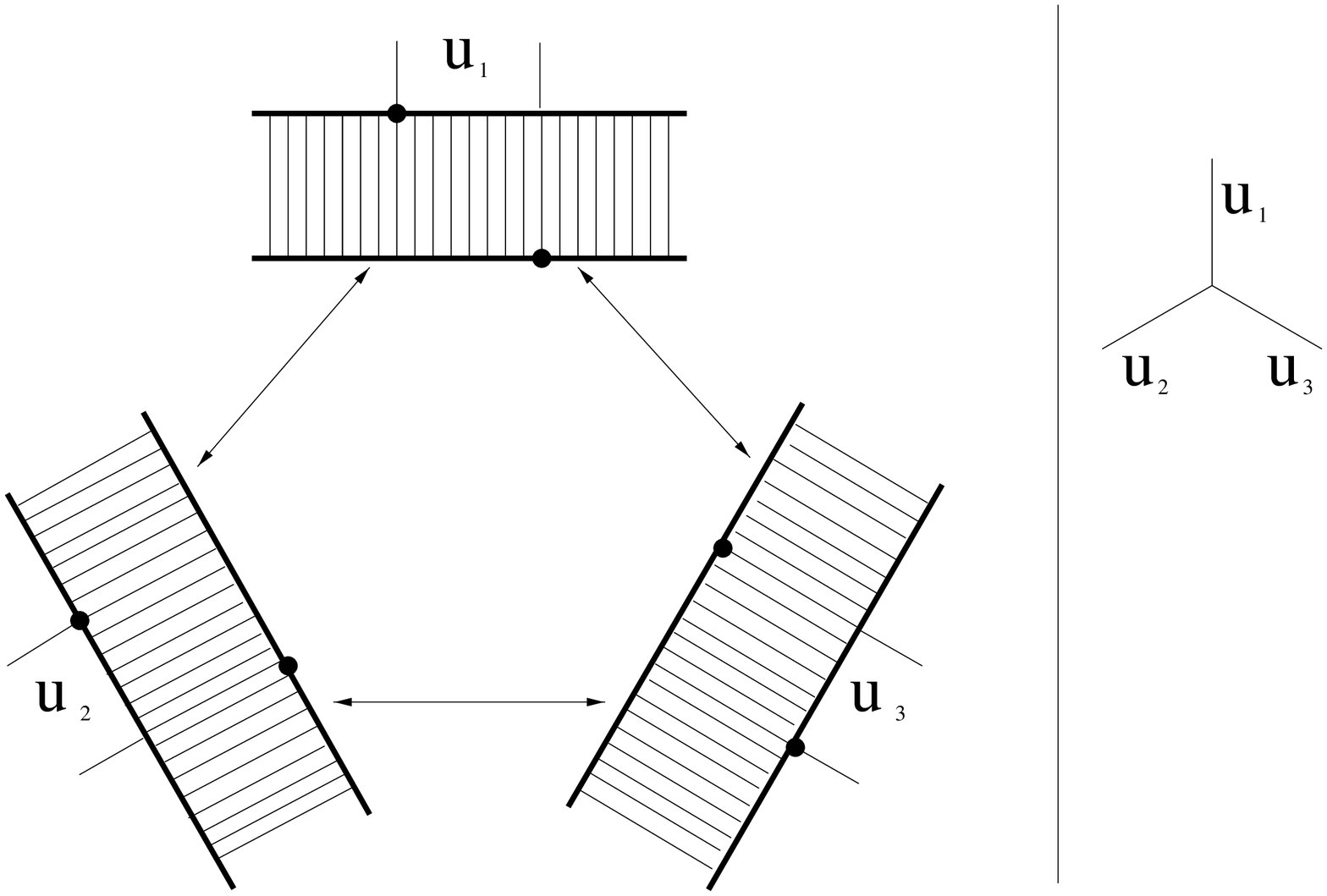}}}
\refstepcounter{equation}\label{ublam-}
\vspace{5mm}
\centerline{Fig. \arabic{equation}}
\vspace{5mm}

 Note that although we have started with infinite bunches of curves the
resulting lamination is finite. All these curves glue together into
a finite number of connected components and possibly infinite number of closed
curves surrounding punctures. Indeed, any curve of the lamination is either
closed or goes diagonally along at least one edge. Since the total number of
pieces of right or left handed curves $I=\vert\sum_{\alpha \in
E(\Gamma)}z_\alpha\vert$ is finite the resulting lamination contain no more
than this number of connected components. (In fact the number of connected
components equals $I$ provided all numbers $z_\alpha$ are all nonpositive or
all nonnegative.)

\paragraph{Graph and orientation changes.}

Here is the transformation law for the constructed coordinates for a flip of an
edge for a simple graph.

\setlength{\unitlength}{0.17mm}%
\vspace{5mm}
\centerline{
\begin{picture}(500,275)(20,485)
\thicklines
\put(280,700){\line( 1,-2){ 40}}
\put(320,620){\line( 1, 0){280}}
\put(600,620){\line( 1, 2){ 40}}
\put(600,620){\line( 1,-2){ 40}}
\put(320,620){\line(-1,-2){ 40}}
\thinlines
\put(180,620){\vector(-1, 0){  0}}
\put(180,620){\vector( 1, 0){ 55}}
\thicklines
\put(100,540){\line( 2,-1){ 80}}
\put(100,700){\line( 0,-1){160}}
\put(100,540){\line(-2,-1){ 80}}
\put( 20,740){\line( 2,-1){ 80}}
\put(100,700){\line( 2, 1){ 80}}
\put( 40,745){\makebox(0,0)[lb]{$a$}}
\put(160,745){\makebox(0,0)[rb]{$b$}}
\put(120,620){\makebox(0,0)[lc]{$z$}}
\put(160,480){\makebox(0,0)[rb]{$c$}}
\put( 40,480){\makebox(0,0)[lb]{$d$}}
\put(300,540){\makebox(0,0)[lt]{$d - \min(z,0)$}}
\put(620,540){\makebox(0,0)[rt]{$c + \max(z,0)$}}
\put(620,690){\makebox(0,0)[rb]{$b - \min(z,0)$}}
\put(300,690){\makebox(0,0)[lb]{$a + \max(z,0)$}}
\put(470,640){\makebox(0,0)[cb]{$-z$}}
\end{picture}}\nopagebreak
\refstepcounter{equation}\label{Hflipl}
\centerline{Fig. \arabic{equation}}
(Only changing part of the graph is shown here, the numbers on the other edges
remain unchanged.)

One can write down explicitly what happens to the coordinates when one changes
the orientation of a hole. Since the formulae are relatively complicated we
postpone them to the sixth section.

\paragraph{Relations and common properties of
${\cal L}^d(S)$ and ${\cal L}^H(S)$.}~

\vspace{3mm}{\bf 1.}  Since the transformation rules for coordinates (\ref{dflipl}) and
(\ref{Hflipl}) are continuous w.r.t.  the standard topology of ${\Bbb Q}^n$ the
coordinates define a natural topology on the lamination spaces.  One now
can define the spaces of {\em real measured laminations} (resp.  bounded and
unbounded) as a completion of the corresponding spaces of rational laminations.
These spaces are denoted as ${\cal L}^d$,${\cal L}^h$ and ${\cal L}^H$, respectively.  
Of course we have the coordinate systems on these spaces automatically.

  Note that to define real measured laminations it is not enough just to
replace rational numbers by real numbers in the definition of the space of
laminations. Such definition would not be equivalent to the one given
above since a sequence of more and more complicated curves with smaller and
smaller weights may converge to a real measured lamination. Thurston in
\cite{Thurston} defined real measured laminations directly as transversely
measured foliation of closed submanifolds. It seems to us that our definition
is more convenient for practical computations although it does not work well
for surfaces without boundary.

\vspace{3mm}{\bf 2.} An unbounded lamination is integral if and only if it has
integral coordinates. A bounded lamination is integral if and only if it
has integral coordinates  and the sum of three numbers on edges incident to
each vertex is even.

\vspace{3mm}{\bf 3.} If $\sigma: \tilde{S} \rightarrow S$ is an unramified covering then
we can define an {\em inverse image} $\sigma^*(f) \in {\cal
L}^d(\tilde{S})$ of a lamination $f \in {\cal L}^d(S)$. For a rational $f$
the curves of $\sigma^*(f)$ are just  full inverse images of the curves of
$f$ with the same weights as on the respective curves of $f$. This mapping
can be obviously extended to all laminations. The analogous mapping
$\sigma^*: {\cal L}^H(S) \rightarrow {\cal L}^H(\tilde{S})$ can be
analogously defined for the spaces of unbounded laminations.

Note that the graph coordinates $\{\tilde{z}_{\tilde{\alpha}}\vert
\tilde{\alpha} \in E(\sigma^*\Gamma)\}$ of a lamination $\sigma^*(f)$ w.r.t.
to the graph $\sigma^*\Gamma$ are just pullbacks of the graph coordinates
$\{z_\alpha\vert \alpha \in E(\Gamma)\}$ of
the lamination $f$ w.r.t.  $\Gamma$, i.e., $\tilde{z}_{\tilde{\alpha}} =
z_{\sigma\tilde{\alpha}}$

  The constructed mappings $\sigma^*(S):{\cal L}^d(S) \rightarrow
{\cal L}^d(\tilde{S})$ and  $\sigma^*(S):{\cal L}^H(S) \rightarrow
{\cal L}^H(\tilde{S})$ are embedings.

\vspace{3mm}{\bf 4.} Denote the closure of ${\Bbb Q}{\cal L}(S)$ in ${\cal L}^d(S)$
by ${\cal L}(S)$. We have the following commuting diagram of natural mappings
commuting with the action of the mapping class group:
\begin{equation}
\setlength{\unitlength}{0.05em}
\begin{picture}(440,155)(190,660)
\thinlines
\put(620,760){\vector( 0,-1){ 80}}
\put(520,780){\vector( 1, 0){70}}
\put(480,760){\vector( 0,-1){ 80}}
\put(500,660){\vector( 1, 0){100}}
\put(360,780){\vector( 1, 0){80}}
\put(130,760){\vector( 0,-1){ 80}}
\put(160,790){\vector( 1, 0){140}}
\put(300,770){\vector(-1, 0){140}}
\put(130,780){\makebox(0,0)[cc]{${\cal L}^d(S)$}}
\put(330,780){\makebox(0,0)[cc]{${\cal L}(S)$}}
\put(480,780){\makebox(0,0)[cc]{${\cal L}^H(S)$}}
\put(620,780){\makebox(0,0)[cc]{${\cal L}^h(S)$}}
\put(620,660){\makebox(0,0)[cc]{${\Bbb R}_+^s$}}
\put(480,660){\makebox(0,0)[cc]{${\Bbb R}^s$}}
\put(130,660){\makebox(0,0)[cc]{${\Bbb R}^s$}}
\put(150,720){\makebox(0,0)[rc]{$a$}}
\put(230,800){\makebox(0,0)[lb]{$p$}}
\put(230,760){\makebox(0,0)[ct]{$i_{v_1,\ldots,v_s}$}}
\put(400,790){\makebox(0,0)[cb]{$i$}}
\put(550,790){\makebox(0,0)[cb]{$\Sigma$}}
\put(550,670){\makebox(0,0)[cb]{$\Sigma_0$}}
\put(470,720){\makebox(0,0)[rc]{$l^H$}}
\put(630,720){\makebox(0,0)[lc]{$l^h$}}
\end{picture}
\end{equation}
The projection $p$ forgets the curves surrounding holes; the projection
$l^h$ (resp $l^H$) is given by the total weights of ends of curves entering the
hole (resp. taken with minus sign for the case of ${\cal L}^H(S)$ if the
orientation of the hole is opposite to the one induced by the orientation of
the surface); $\Sigma$ and $\Sigma_0$ are the canonical projections on the
quotient by the group $({\Bbb Z}/2{\Bbb Z})^s$  acting by changing orientation
 of the holes on
${\cal L}^H(S)$ and by changing sign of the standard coordinates on ${\Bbb R}^s$,
respectively; $a$ is given by the weights of the curves surrounding holes;
$i(v_1,\ldots,v_s)$ is a family of embeddings characterised by the condition
that $ai_{v_1,\ldots,v_s}(x) =  (v_1,\ldots,v_s )$ for any $x \in {\cal L}(S)$.
The image of $i$ coincides with the kernel of $l^H$ and with the stable points of
the $({\Bbb Z}/2{\Bbb Z})^s$ action.

In coordinates the mapping $ip : {\cal L}^d(S) \rightarrow {\cal L}^H(S)$ is
given by
\begin{equation}
z_\aalpha = u_{\aalpha(1)} + u_{\aalpha(3)} - u_{\aalpha(2)} - u_{\aalpha(4)}
\end{equation}
where $\{z_\alpha | \alpha \in E(\gamma)\}$ and $\{u_\alpha| \alpha \in
E(\gamma)\}$ are the coordinates on ${\cal L}^H(S)$ and ${\cal L}^d(S)$,
respectively, w.r.t. the same graph $\Gamma$. By $z_\aalpha$  and $u_\aalpha$ we 
mean the numbers assigned to the corresponding unoriented edges.

The mapping $a$ is given by
\begin{equation}
\{u_\alpha| \alpha \in EE(\gamma)\} \mapsto \{\frac{1}{2}\max_{\aalpha \in
\gamma} (-u_\aalpha + u_{\aalpha(1)} - u_{\aalpha(4)})| \gamma \in
F(\Gamma)\}\label{col} \end{equation}

The mapping $l^H$ (resp. $l^h$) is given by
\begin{equation}
\{z_\alpha| \alpha \in E(\gamma)\} \mapsto \{\sum_{\aalpha \in \gamma} z_\aalpha
| \gamma \in F(\Gamma)\} {\rm ~~(resp.,~~}\{z_\aalpha| \aalpha \in E(\gamma)\}
\mapsto \{\vert \sum_{\aalpha \in \gamma} z_\aalpha\vert | \gamma \in
F(\Gamma)\} {\rm )}
\end{equation}

\section{Teichm\" uller spaces.}

{\em The Teichm\" uller space} ${\cal T}(S)$ (resp. {\em Moduli space} ${\cal
M}(S)$) of a closed surface $S$ is the space of complex structures on $S$ modulo
diffeomorphisms homotopy equivalent to the identity (resp. modulo all
diffeomorphisms). The extension of these notions to surfaces with boundary
depends on the condition that one imposes on the behaviour of the
complex structure at the boundary of the surface. The most traditional
definition considers only complex structures degenerating at the boundary and
such that a tubular neighbourhood of each boundary component is isomorphic as a
complex manifold to a punctured disc. (Such kind of singularity is called {\em
puncture}.) We denote the corresponding Teichm\" uller and moduli spaces by the
same letters ${\cal T}(S)$ and ${\cal M}(S)$, respectively. We describe two
other modifications of Teichm\" uller spaces and give explicit
parameterisations of them.
But before we just recall some basic facts about relations between complex
structures, constant negative curvature metrics and discrete subgroups of the
group $PSL(2,{\Bbb R})$. For more details we recommend the reviews \cite{KAG}.

 According to the Poincar\' e uniformisation theorem any complex surface $S$
can be represented as a quotient of the upper half plane $H$ by a discrete
subgroup $\Delta$ of its automorphism group (sometimes called the {M\" obius
group}) $PSL(2,{\Bbb R})$ of real $2 \times 2$ matrices with unit
determinant considered up to the factor $ -1$.  The group $\Delta$ is
canonically isomorphic to the fundamental group of the surface $\pi_1(S)$.
$\Delta$ is defined by the complex structure of the surface up to conjugation
by an element of $PSL(2,{\Bbb R})$.  Therefore we get an embedding ${\cal
T}(S) \mapsto {\rm Hom}(\pi_1(S) \rightarrow PSL(2,{\Bbb R}))/PSL(2,{\Bbb
R})$.  This embedding has the following properties:

1. The image of any loop is a matrix with one or two real eigenvectors. (Such
elements of $PSL(2,{\Bbb R})$ are called {\em parabolic} and {\em hyperbolic},
respectively.)

2. Parabolic elements correspond to loops surrounding punctures only.

The proof of these well known properties will in particular follow from the
construction of the parameterisations.

On $H$ there exists a unique $PSL(2,{\Bbb R})$-invariant Riemann curvature $-1
$ metric. It induces a metric on $S$. Since this metric is of negative
curvature, any homotopy class of closed curves contains a unique geodesics.
Homotopy classes of closed curves on a surface are in one-to-one correspondence
with the conjugacy classes of its fundamental group $\Delta$. Denote by
$\gamma$ an element of $\pi_1(S)$ and by $l(\gamma)$ the length of the
corresponding geodesics. Then a simple computation shows, that

\begin{equation}
l(\gamma) = \left|\log\frac{\lambda_1}{\lambda_2}\right|,\label{length}
\end{equation}
where $\lambda_1$ and $\lambda_2$ are the eigenvalues of the element of
$PSL(2,{\Bbb R})$ corresponding to $\gamma$. This number is obviously
correctly
defined, i.e., it does not depend on the choices of particular representation
of $\pi_1(S)$, of a particular element of $\pi_1(S)$, representing a given loop
and of a particular $2\times 2$ matrix representing an element of
$PSL(2,{\Bbb R})$. This formula implies that the length of a
geodesics surrounding a puncture is zero. Note that taking curvature to be $-1$ is
equivalent to the demand that the curvature is negative and constant and areas
of  ideal triangles are equal to $\pi$, which
normalisation condition is more convenient practically.

\subsection{Teichm\" uller space of surfaces with holes ${\cal T}^H(S)$.}

\paragraph{Definition.}
There is another condition one can impose on the behaviour of complex structure
in a vicinity of the boundary and still get a finite dimensional moduli
space. Demand that a boundary component be either a puncture or the
complex structure is nondegenerate at the boundary. A boundary of the latter
type is called a {\em hole.} A neighbourhood of a hole is isomorphic as a
complex manifold to an annulus.  The corresponding moduli space is denoted by
${\cal T}^h(S)$.

For our purposes it is more convenient to introduce another space. ${\cal
T}^H(S)$ is the space of complex structures on $S$ together with orientations of
all holes. (By orientation of a hole we mean the orientation of
the corresponding boundary component.) Note that, although it is not {\em a
priori} obvious, this space possesses a natural topology in which it is connected.

\paragraph{Construction of coordinates.}

Let $\Gamma \in \Gamma^3(S)$ be a three-valent graph,
corresponding to a surface $S$. For any point of ${\cal T}^H(S)$ we are going
to describe a rule for assigning a real number to each edge of $\Gamma$. The
collection of these numbers  will give us a global parameterisation of ${\cal
T}^H(S)$.

For simplicity consider first the case, when all boundary components are
holes. Draw a closed geodesics around each hole and  cut out cylinders by
them. We thus get a surface with geodesic boundary. Then cut the surface by
the edges of the dual graph $\Gamma^\vee$ into hexagons. (These edges are not
necessarily geodesic though one can suppose them to be.)  Take an edge and
two hexagons incident to it and lift the resulting octagon to the upper half
plane $H$. The octagon has four geodesic sides facing
holes. Continue these geodesics up to the real axis. Now, the orientations
of the holes induce the orientations of the geodesics. Using these
orientations choose one of the two infinities of each geodesics, say, the end.
We obtain therefore four
points on the real axis, or to be more precise, on ${\Bbb R}P^1$. Note, that
the four geodesics do not intersect on $H$, and therefore the cyclic order of
the constructed points on ${\Bbb R}P^1$ does not depend on the point of ${\cal
T}^H(S)$ we have started with. Among the constructed four points there are two
distinguished ones which originate from the geodesics connected by the
edge we have started with. Using the action of the M\" obius group on the upper half
plane, we can shift these two points to zero and infinity, respectively, and one
of the remaining points to $-1$. And now finally assign to the edge the
logarithm of the coordinate of the fourth point. (Of course this forth
coordinate is nothing but a suitable cross-ratio of those four points.)

\vspace{5mm}
\centerline{\epsfxsize30\baselineskip\epsfbox{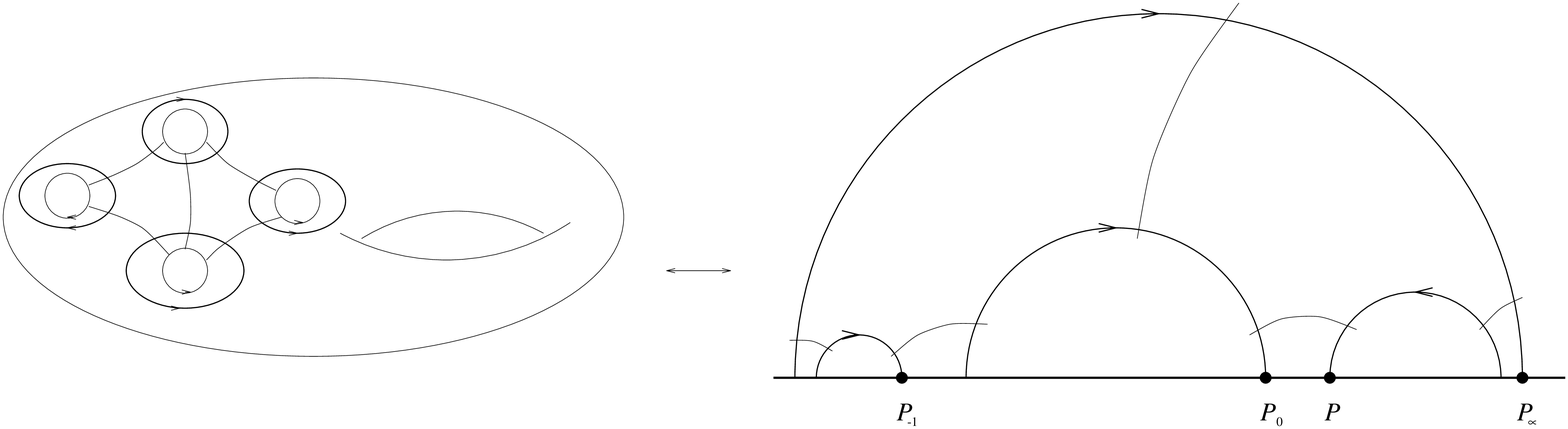}} 
\refstepcounter{equation}\label{Hcoord}
\vspace{5mm}
\centerline{Fig. \arabic{equation}}
\vspace{5mm}

Note that if we have punctures instead of some holes it does not spoil the
construction.  In this case some edges of the considered hexagons shrink to
points, the corresponding geodesics on the upper half plane shrink to points
on the real axis and no orientation is necessary to choose between their ends.

\paragraph{Reconstruction.}
Our goal now is to construct a surface starting from a three-valent fat graph
$\Gamma \in \Gamma^3(S)$ with real numbers $\{z_\alpha\vert \alpha \in
E(\Gamma)\}$ on edges.  First of all give a simple receipt how to restore
orientations of the boundary components from these data: The orientation of a
boundary component corresponding to a face $\gamma$ is just induced from the
orientation of the surface (resp. opposite to the induced one) if the sum
$\sum_{\aalpha \in \gamma}z_\aalpha$ is positive (resp. negative). If the sum
is zero, it means, that it is not a boundary, but a puncture.

Construction of the surface itself can be achieved in two equivalent ways. We
shall describe both since one is more transparent from the geometric point of view
and the other is useful for practical computations.

{\em Construction by gluing.} We are going to glue our surface out of ideal
hyperbolic triangles. The lengths of the sides of ideal triangles are
infinite and therefore we can glue two triangles in many ways which differ by
shifting one triangle w.r.t. another along the side. The ways of gluing
triangles can be parameterised by the cross-ratios of four vertices of the
obtained quadrilateral (considered as points of ${\Bbb R}P^1$). For our
purpose it is convenient to take as a parameter $z$ the logarithm of the
cross-ratio
\begin{equation}
z = \log\frac{(P_0-P)(P_{-1}-P_\infty)}{(P_\infty - P)(P_{-1}-P_0)},
\end{equation}
where $P_0, P_{-1}, P$, and $P_\infty$ are coordinates of vertices of the
quadrilateral, $P_0$ and $P_\infty$ being coordinates of the ends of the side
we are gluing triangles along.

Now consider the dual graph $\Gamma^\vee$. Its faces are triangles. Take one
ideal hyperbolic triangle for each face of this graph and glue them together
along the edges just as they are glued in $\Gamma^\vee$ using numbers assigned
to the edges as gluing parameters.

Note that although this is not quite obvious the resulting surface is not
necessarily complete. In fact it is not the original surface with the absolute
boundary but only what we get out of it by cutting off annuli around holes by
closed geodesics.

{\em Construction of the Fuchsian group}. We are now going to construct a
discrete subgroup $\Delta$ of $PSL(2,{\Bbb R})$ starting from a graph
$\Gamma \in \Gamma^3(S)$ with numbers on edges. Modify first the original graph $\Gamma$ at
each vertex in the following way. Disconnect the  edges at the vertex
and then connect them by three more edges forming a triangle. Orient
the edges of the triangle in the counterclockwise direction. Now assign to each
of these edges the matrix
$I = \left(\begin{array}{cc} 1 & 1 \\ -1 & 0 \end{array}
\right)$.
Assign to each old edge $\alpha$ the matrix
$A(z_\alpha) = \left(\begin{array}{cc} 0 &
e^{z_\alpha/2} \\ -e^{-z_\alpha/2} & 0 \end{array} \right)$.
Now for any oriented path on this graph we can associate a matrix by
multiplying consecutively all matrices we meet along it, taking $I^{-1}$
instead of $I$ each time we go along a new edge in inverse direction w.r.t.
the orientation. (The orientation of the old edges is not to be taken into
account, since $A(z)^2=-1$ and therefore $A(z)$ coincides with its inverse in
the group $PSL(2,{\Bbb R})$.) In particular, if we take closed paths starting
form a fixed vertex of the graph, we get a homomorphism of the fundamental
group of $\Gamma$ to the group $PSL(2,{\Bbb R})$. The image of this
homomorphism is just the desired group $\Delta$.

 In principle we need to prove that these two constructions are inverse to
the above construction of coordinates indeed, what is almost obvious,
especially for the first one.  The only note we would like to make here is to
show where the matrices $I$ and $A(z)$ came from.  Consider two ideal
triangles on the upper half plane with vertices at the points ${-1,0,\infty}$
and ${e^{z},\infty,0}$, respectively. Then the M\" obius transform which
permutes the vertices of the first triangle is given by the matrix $I$, and
the one which maps one triangle to another (respecting the order of vertices
given two lines above) is given by $A(z)$.

\paragraph{Graph and orientation change.}
Here is the transformation law for the constructed coordinates for a flip of an
edge.

\centerline{
\setlength{\unitlength}{0.17mm}%
\begin{picture}(500,275)(20,485)
\thicklines
\put(280,700){\line( 1,-2){ 40}}
\put(320,620){\line( 1, 0){300}}
\put(620,620){\line( 1, 2){ 40}}
\put(620,620){\line( 1,-2){ 40}}
\put(320,620){\line(-1,-2){ 40}}
\thinlines
\put(180,620){\vector(-1, 0){  0}}
\put(180,620){\vector( 1, 0){ 55}}
\thicklines
\put(100,540){\line( 2,-1){ 80}}
\put(100,700){\line( 0,-1){160}}
\put(100,540){\line(-2,-1){ 80}}
\put( 20,740){\line( 2,-1){ 80}}
\put(100,700){\line( 2, 1){ 80}}
\put( 40,745){\makebox(0,0)[lb]{$a$}}
\put(160,745){\makebox(0,0)[rb]{$b$}}
\put(120,620){\makebox(0,0)[lc]{$z$}}
\put(160,485){\makebox(0,0)[rt]{$c$}}
\put( 40,485){\makebox(0,0)[lt]{$d$}}
\put(300,540){\makebox(0,0)[lt]{$d - \log(e^{-z}+1)$}}
\put(640,540){\makebox(0,0)[rt]{$c + \log(e^z+1)$}}
\put(640,690){\makebox(0,0)[rb]{$b - \log(e^{-z}+1)$}}
\put(300,690){\makebox(0,0)[lb]{$a + \log(e^z+1)$}}
\put(480,640){\makebox(0,0)[cb]{$-z$}}
\label{Hflipt}
\end{picture}
}
\vspace{5mm}
\centerline{Fig. \arabic{equation}}
(Only changing part of the graph is shown here, the numbers on the other edges
remain unchanged.)

The formulae for changing orientations of holes are slightly more complicated.
We reproduce here how the coordinates transform under change of orientation of
one hole $\gamma$ on a simple graph.

\begin{equation}
z_{\alpha_0} \rightarrow - \frac{1}{2}(z_{\alpha_1} +
z_{\alpha_n}) - [z_{\alpha_2},\ldots,z_{\alpha_n},z_{\alpha_0}] +
[z_{\alpha_0},\ldots,z_{\alpha_{n-1}}],
 \label{ori1}
\end{equation}
where $\alpha_0,\ldots,\alpha_n$ are the edges belonging to $\gamma$
numerated in counterclockwise order.

If $\aalpha \notin \gamma$, but an end of $\alpha$ belongs to $\gamma$ then
the number on it also changes
\begin{equation} z_\alpha \rightarrow z_\alpha
+ \frac{1}{2}(z_{\alpha_0} + z_{\alpha_n}) +
[z_{\alpha_1},\ldots,z_{\alpha_n}] - [z_{\alpha_0},\ldots,z_{\alpha_{n-1}}],
\label{ori2}
\end{equation}
where $\alpha_0,\ldots,\alpha_n$ are the edges belonging to $\gamma$
numerated in counterclockwise order starting with the intersection point with
$\alpha$.

We have used the notation
\begin{equation}
\label{inv1} [x_1,\ldots,x_i] := \log (e^{(x_1+\cdots+x_n)/2} +
e^{(x_1+\cdots+x_{n-1}-x_n)/2} + \ldots + e^{(-x_1 - \cdots - x_n)/2})
\end{equation}

The promised formulae describing how graph coordinates of an unbounded
lamination change under the change of orientation of a hole can be obtained
from (\ref{ori1})---(\ref{ori2}) by taking the limit $z_\alpha \rightarrow
\infty$ or in another words by replacing (\ref{inv1}) by
\begin{equation}
[x_1,\ldots,x_i] := \frac{1}{2}\max({(x_1+\cdots+x_n)},
{(x_1+\cdots+x_{n-1}-x_n)} , \ldots , {(-x_1 - \cdots - x_n)})
\end{equation}

If $\alpha$ does not intersect with $\gamma$ then $z_\alpha$ does not change

For a nonregular graph the formulae are slightly more complicated. For example
for a torus with one hole they are

\begin{equation}
x \rightarrow -x + 2[z,x,0,x,y,0,y,z]-2[y,z,0,z,x,0,y]
\end{equation}
and analogously for the other two edges.

\subsection{Decorated Teichm\" uller space ${\cal T}^d(S)$.}
 In this section we are going to reproduce some results of Penner \cite{Penner}.

 Before giving a definition of the decorated Teichm\" uller space ${\cal
T}^d(S)$, recall what a horocycle around a puncture on a Riemann surface is.
Consider the upper half plane $H$ with the standard curvature $-1$ metric.
Then a horocycle is either a circle tangent to the real axis or a horizontal
line.  Another more scientific definition is the following one. Let $O(x,y)$,
where $x,y \in H$, be the set of points equidistant from $x$ and containing
$y$. The limit of $O(x,y)$ when $x$ tends to a point $p$ at infinity is
called a horocycle and the point $p$ is called its {\em basepoint}. The space
of horocycles based at a given point is homeomorphic to a real line. A
horocycle based at a point $p$ is setwise stable under the action of the
subgroup of parabolic elements of $PSL(2,{\Bbb R})$ stabilising $p$.

Now consider a punctured Riemann surface $S = H/\Delta$, where $\Delta$ is a
discrete subgroup of $PSL(2,{\Bbb R})$. Consider a point $p$ stabilised by a
parabolic element of $\Delta$ and a horocycle based at $p$. A horocycle on $S$
by definition is the image of such horocycle on $H$. If the original horocycle
on $H$ is small enough, its image on $S$ is a small circle surrounding a
puncture and orthogonal to any geodesics coming out of this puncture. But a
projection general horocycle to the surface may have a relatively complicated topology.

\paragraph{Definition.}
  A {\em decorated Riemann surface} is a punctured Riemann surface with a
horocycle chosen around each puncture. The Teichm\" uller space of decorated
surfaces is called the {\em decorated Teichm\" uller space} and is denoted by
${\cal T}^d(S)$.

\paragraph{Construction of coordinates.}
Take a three-valent graph, corresponding to the surface $\Gamma \in
\Gamma_0^3(S)$. Assign now a real number to each edge of $\Gamma$. Take an
edge of the dual graph $\Gamma^\vee$ corresponding to a given edge of $\Gamma$.
This edge connects two punctures of $S$. Make it geodesic by a homotopy. Now consider the inverse image of this geodesics together with the horocycles around its
endpoints on the upper half plane. Assign now to the edge of the original
graph the length of the part of the geodesics on $H$ between two horocycles if
the latter ones do not intersect. If they do, assign the length with the minus
sign. (fig. \ref{copenner})

\centerline{ 
\unitlength 0.8mm
\linethickness{0.4pt}
\begin{picture}(169.33,50.67)(0,80)
\multiput(20.00,95.00)(-0.15,0.12){13}{\line(-1,0){0.15}}
\multiput(18.09,96.51)(-0.13,0.12){13}{\line(-1,0){0.13}}
\multiput(16.38,98.01)(-0.12,0.12){13}{\line(0,1){0.12}}
\multiput(14.88,99.52)(-0.12,0.14){11}{\line(0,1){0.14}}
\multiput(13.58,101.03)(-0.11,0.15){10}{\line(0,1){0.15}}
\multiput(12.48,102.54)(-0.11,0.19){8}{\line(0,1){0.19}}
\multiput(11.58,104.04)(-0.12,0.25){6}{\line(0,1){0.25}}
\multiput(10.88,105.55)(-0.10,0.30){5}{\line(0,1){0.30}}
\multiput(10.38,107.06)(-0.10,0.50){3}{\line(0,1){0.50}}
\put(10.09,108.57){\line(0,1){1.51}}
\put(10.00,110.07){\line(0,1){1.51}}
\multiput(10.11,111.58)(0.10,0.50){3}{\line(0,1){0.50}}
\multiput(10.42,113.09)(0.10,0.30){5}{\line(0,1){0.30}}
\multiput(10.94,114.60)(0.12,0.25){6}{\line(0,1){0.25}}
\multiput(11.66,116.10)(0.11,0.19){8}{\line(0,1){0.19}}
\multiput(12.58,117.61)(0.11,0.15){10}{\line(0,1){0.15}}
\multiput(13.70,119.12)(0.11,0.13){12}{\line(0,1){0.13}}
\multiput(15.02,120.63)(0.12,0.12){13}{\line(1,0){0.12}}
\multiput(16.54,122.13)(0.13,0.12){13}{\line(1,0){0.13}}
\multiput(18.27,123.64)(0.14,0.11){12}{\line(1,0){0.14}}
\multiput(20.00,95.00)(0.28,-0.12){8}{\line(1,0){0.28}}
\multiput(22.27,94.07)(0.32,-0.12){7}{\line(1,0){0.32}}
\multiput(24.55,93.25)(0.38,-0.12){6}{\line(1,0){0.38}}
\multiput(26.82,92.55)(0.45,-0.12){5}{\line(1,0){0.45}}
\multiput(29.09,91.95)(0.45,-0.10){5}{\line(1,0){0.45}}
\multiput(31.36,91.47)(0.57,-0.09){4}{\line(1,0){0.57}}
\multiput(33.64,91.10)(0.76,-0.09){3}{\line(1,0){0.76}}
\multiput(35.91,90.85)(1.14,-0.07){2}{\line(1,0){1.14}}
\put(38.18,90.70){\line(1,0){2.27}}
\put(40.45,90.67){\line(1,0){2.27}}
\multiput(42.73,90.75)(1.14,0.10){2}{\line(1,0){1.14}}
\multiput(45.00,90.94)(0.76,0.10){3}{\line(1,0){0.76}}
\multiput(47.27,91.24)(0.57,0.10){4}{\line(1,0){0.57}}
\multiput(49.55,91.65)(0.45,0.11){5}{\line(1,0){0.45}}
\multiput(51.82,92.18)(0.38,0.11){6}{\line(1,0){0.38}}
\multiput(54.09,92.82)(0.32,0.11){7}{\line(1,0){0.32}}
\multiput(56.36,93.57)(0.30,0.12){12}{\line(1,0){0.30}}
\multiput(60.00,125.00)(0.15,-0.12){13}{\line(1,0){0.15}}
\multiput(61.91,123.49)(0.13,-0.12){13}{\line(1,0){0.13}}
\multiput(63.62,121.99)(0.12,-0.12){13}{\line(0,-1){0.12}}
\multiput(65.12,120.48)(0.12,-0.14){11}{\line(0,-1){0.14}}
\multiput(66.42,118.97)(0.11,-0.15){10}{\line(0,-1){0.15}}
\multiput(67.52,117.46)(0.11,-0.19){8}{\line(0,-1){0.19}}
\multiput(68.42,115.96)(0.12,-0.25){6}{\line(0,-1){0.25}}
\multiput(69.12,114.45)(0.10,-0.30){5}{\line(0,-1){0.30}}
\multiput(69.62,112.94)(0.10,-0.50){3}{\line(0,-1){0.50}}
\put(69.91,111.43){\line(0,-1){1.51}}
\put(70.00,109.93){\line(0,-1){1.51}}
\multiput(69.89,108.42)(-0.10,-0.50){3}{\line(0,-1){0.50}}
\multiput(69.58,106.91)(-0.10,-0.30){5}{\line(0,-1){0.30}}
\multiput(69.06,105.40)(-0.12,-0.25){6}{\line(0,-1){0.25}}
\multiput(68.34,103.90)(-0.11,-0.19){8}{\line(0,-1){0.19}}
\multiput(67.42,102.39)(-0.11,-0.15){10}{\line(0,-1){0.15}}
\multiput(66.30,100.88)(-0.11,-0.13){12}{\line(0,-1){0.13}}
\multiput(64.98,99.37)(-0.12,-0.12){13}{\line(-1,0){0.12}}
\multiput(63.46,97.87)(-0.13,-0.12){13}{\line(-1,0){0.13}}
\multiput(61.73,96.36)(-0.14,-0.11){12}{\line(-1,0){0.14}}
\bezier{176}(60.00,125.00)(40.00,133.67)(20.00,125.00)
\bezier{132}(35.00,110.00)(47.33,99.33)(60.00,110.00)
\bezier{72}(40.00,106.33)(47.33,111.00)(55.00,106.33)
\put(25.00,120.00){\circle{10.00}}
\put(25.00,100.00){\circle{10.00}}
\put(25.00,100.00){\circle*{1.33}}
\put(25.00,120.00){\circle*{1.33}}
\put(25.00,120.00){\line(0,-1){20.00}}
\bezier{200}(109.33,95.00)(89.33,110.00)(109.33,125.00)
\bezier{176}(109.33,95.00)(129.33,86.33)(149.33,95.00)
\bezier{200}(149.33,125.00)(169.33,110.00)(149.33,95.00)
\bezier{176}(149.33,125.00)(129.33,133.67)(109.33,125.00)
\bezier{132}(124.33,110.00)(136.67,99.33)(149.33,110.00)
\bezier{72}(129.33,106.33)(136.67,111.00)(144.33,106.33)
\put(114.33,113.00){\circle{10.00}}
\put(114.33,106.00){\circle{10.00}}
\put(114.33,106.00){\circle*{1.33}}
\put(114.33,113.00){\circle*{1.33}}
\put(114.33,106.00){\line(0,1){7.00}}
\put(114.33,113.00){\line(0,1){0.00}}
\put(40.00,80.00){\makebox(0,0)[cc]{$+$}}
\put(130.00,80.00){\makebox(0,0)[cc]{$-$}}
\end{picture}
}
\vspace{5mm}
\refstepcounter{equation}\label{copenner}
\centerline{Fig. \arabic{equation}}

\paragraph{Reconstruction} is quite analogous to that for holed
surfaces. There is a canonical mapping $ip:{\cal T}^d(S) \rightarrow {\cal
T}^H(S)$ which just forgets about horocycles and will be given explicitly in
coordinates by (\ref{tembed}).
Therefore, to reconstruct the surface itself we can just apply the
reconstruction procedure for ${\cal T}^H(S)$. To reconstruct the horocycles
consider an ideal triangle which we have used to glue the surface. On each
edge we have a length of the corresponding geodesics between the horocycles. It
allows us to restore unambiguously the points of intersection of the horocycles
with the edges.

\paragraph{Graph change.}
\centerline{
\setlength{\unitlength}{0.17mm}%
\begin{picture}(500,275)(20,485)
\thicklines
\put(280,700){\line( 1,-2){ 40}}
\put(320,620){\line( 1, 0){280}}
\put(600,620){\line( 1, 2){ 40}}
\put(600,620){\line( 1,-2){ 40}}
\put(320,620){\line(-1,-2){ 40}}
\thinlines
\put(180,620){\vector(-1, 0){  0}}
\put(180,620){\vector( 1, 0){ 55}}
\thicklines
\put(100,540){\line( 2,-1){ 80}}
\put(100,700){\line( 0,-1){160}}
\put(100,540){\line(-2,-1){ 80}}
\put( 20,740){\line( 2,-1){ 80}}
\put(100,700){\line( 2, 1){ 80}}
\put( 40,745){\makebox(0,0)[lb]{$A$}}
\put(160,745){\makebox(0,0)[rb]{$B$}}
\put(120,620){\makebox(0,0)[lc]{$Z$}}
\put(160,485){\makebox(0,0)[rt]{$C$}}
\put( 40,485){\makebox(0,0)[lt]{$D$}}
\put(300,540){\makebox(0,0)[lt]{$D$}}
\put(620,540){\makebox(0,0)[rt]{$C$}}
\put(620,690){\makebox(0,0)[rb]{$B$}}
\put(300,690){\makebox(0,0)[lb]{$A$}}
\put(470,640){\makebox(0,0)[cb]{$\log(e^{A+C}+e^{B+D})-Z$}}
\end{picture}
}
\refstepcounter{equation}\label{dflipt}
\vspace{5mm}
\centerline{Fig. \arabic{equation}}
\vspace{5mm}
(Only changing part of the graph is shown here, the numbers on the other edges
remain unchanged.)

\paragraph{Relations and common properties of ${\cal T}^d$ and ${\cal T}^H$.}~

\vspace{3mm}{\bf 1.}
If $\sigma: \tilde{S} \rightarrow S$ is an unramified covering then
we can define an {\em inverse image} $\sigma^*(m) \in {\cal
T}^d(\tilde{S})$ of a complex structure $m \in {\cal L}^d(S)$ as the unique
complex structure on $\tilde{S}$ s.t. $\sigma$ is holomorphic. Such
mapping can be obviously extended to all laminations. The analogous mapping
$\sigma^*:  {\cal L}^H(S) \rightarrow {\cal L}^H(\tilde{S})$ can be
analogously defined for the spaces of unbounded laminations.

Note that the graph coordinates $\{\tilde{z}_{\tilde{\alpha}}\vert
\tilde{\alpha} \in E(\tilde{\Gamma})\}$ of a complex structure $\sigma^*(m)$
w.r.t.  to the graph $\tilde{\Gamma}$ are just pullbacks of the graph
coordinates $\{z_\alpha\vert \alpha \in E(\Gamma)\}$ of the complex
structure $m$ w.r.t.  $\Gamma$, i.e., $\tilde{z}_{\tilde{\alpha}} =
z_{\sigma\tilde{\alpha}}$

  The constructed mappings $\sigma^*(S):{\cal T}^d(S) \rightarrow {\cal
T}^d(\tilde{S})$ and  $\sigma^*(S):{\cal T}^H(S) \rightarrow {\cal
T}^H(\tilde{S})$ are obviously embedings.

\vspace{3mm}{\bf 2.}
There exists (as for laminations) a set of morphisms between
different versions of Teichm\" uller spaces commuting with the action of the
mapping class group and satisfying analogous properties. We have the
following commuting diagram of natural mappings commuting with the action of
the mapping class group (for simplicity we denote mappings between Teichm\"
uller spaces by the same letter as for the corresponding mappings of the
lamination spaces):
\begin{equation} \setlength{\unitlength}{0.05em}
\begin{picture}(440,155)(190,660)
\thinlines
\put(620,760){\vector( 0,-1){ 80}}
\put(520,780){\vector( 1, 0){70}}
\put(480,760){\vector( 0,-1){ 80}}
\put(500,660){\vector( 1, 0){100}}
\put(360,780){\vector( 1, 0){80}}
\put(130,760){\vector( 0,-1){ 80}}
\put(160,790){\vector( 1, 0){140}}
\put(300,770){\vector(-1, 0){140}}
\put(130,780){\makebox(0,0)[cc]{${\cal T}^d(S)$}}
\put(330,780){\makebox(0,0)[cc]{${\cal T}(S)$}}
\put(480,780){\makebox(0,0)[cc]{${\cal T}^H(S)$}}
\put(620,780){\makebox(0,0)[cc]{${\cal T}^h(S)$}}
\put(620,660){\makebox(0,0)[cc]{${\Bbb R}_+^s$}}
\put(480,660){\makebox(0,0)[cc]{${\Bbb R}^s$}}
\put(130,660){\makebox(0,0)[cc]{${\Bbb R}^s$}}
\put(150,720){\makebox(0,0)[rc]{$a$}}
\put(230,800){\makebox(0,0)[lb]{$p$}}
\put(230,760){\makebox(0,0)[ct]{$i_{v_1,\ldots,v_s}$}}
\put(400,790){\makebox(0,0)[cb]{$i$}}
\put(550,790){\makebox(0,0)[cb]{$\Sigma$}}
\put(550,670){\makebox(0,0)[cb]{$\Sigma_0$}}
\put(470,720){\makebox(0,0)[rc]{$l^H$}}
\put(630,720){\makebox(0,0)[lc]{$l^h$}}
\end{picture}
\end{equation}
The projection $p$ forgets the horocycles; the projection $l^h$
(resp. $l_H$) is given by the lengths of geodesics surrounding the holes (resp.
taken with the minus sign if the orientation of the hole is opposite to the one
induced by the orientation of the surface); $\Sigma$ and $\Sigma_0$ are the
canonical projection on the quotient by the group  acting by changing of
orientations of the holes on ${\cal T}^H(S)$ and by changing signs of
coordinates on ${\Bbb R}^s$, respectively; $a$ is given by the logarithms of
areas of punctured disks bounded by horocycles; $i(v_1,\ldots,v_s)$ is a family
of embeddings characterised by the condition that $ai_{v_1,\ldots,v_s}(x) =
(v_1,\ldots,v_s )$ for any $x \in {\cal T}(S)$. The image of $i$ coincides with
the kernel of $l^H$ and with the stable points of the $({\Bbb Z}/2{\Bbb Z})^s$
action.

In coordinates the mapping $ip : {\cal T}^d(S) \rightarrow {\cal T}^H(S)$ is
given by
\begin{equation}
z_\aalpha = u_{\aalpha(1)} + u_{\aalpha(3)} - u_{\aalpha(2)} - u_{\aalpha(4)}
\label{tembed}
\end{equation}
where $\{z_\alpha | \alpha \in E(\Gamma)\}$ and $\{u_\alpha| \alpha \in
E(\gamma)\}$ are the coordinates on ${\cal T}^H(S)$ and ${\cal T}^d(S)$,
respectively, w.r.t. the same graph $\Gamma$.

The mapping $a$ is given by
\begin{equation}
\{u_\alpha| \alpha \in E(\Gamma)\} \mapsto \{\frac{1}{2}\log\sum_{\aalpha \in \gamma}
e^{-u_\aalpha + u_{\aalpha(1)} - u_{\aalpha(4)}}| \gamma \in
F(\Gamma)\}\label{area}
\end{equation}

The mapping $l^H$ (resp. $l^h$) is given by
\begin{equation}
\{z_\alpha| \alpha \in E(\gamma)\} \mapsto \{\sum_{\alpha \in \gamma} z_\alpha
| \gamma \in F(\Gamma)\} {\rm ~~(resp.~~}\{z_\alpha| \alpha \in E(\gamma)\}
\mapsto \{\left| \sum_{\alpha \in \gamma} z_\alpha\right| | \gamma \in
F(\Gamma)\} {\rm )}
\end{equation}

\paragraph{Relations between Teichm\" uller and lamination spaces.}~

The rules (\ref{dflipt}) and (\ref{dflipl}) show that although the
coordinatewise identification of the spaces ${\cal T}^d(S)$ and ${\cal
L}^d(S)$ is not canonical, i.e., depends on the particular graph, this
identification is canonical asymptotically for large values of coordinates.
In particular it shows, that a projective compactification of the space
${\cal L}^d(S)$ can serve as a compactification boundary for the space ${\cal
T}^d(S)$. This compactification is called {\em Thurston compactification} of
the Teichm\"uller space ${\cal T}^d(S)$. Analogously the rules (\ref{Hflipt})
and  (\ref{Hflipl}) show that the spaces ${\cal L}^H(S)$ and ${\cal T}^H(S)$
are asymptotically canonically isomorphic.
  Further relations will be explained below.

\section{Length of a lamination.}
Suppose we have both complex structure $m \in {\cal T}^h(S)$ and a
rational bounded measured lamination $f \in {\Bbb Q}{\cal L}^d(S)$ on a surface. For each
complex structure we  can associate the constant curvature $-1$ metric on
the surface. We can  deform each curve of the lamination to make it geodesic
and then take a weighted sum of their lengths. This procedure defines a function
$l_{{\cal T}{\cal L}}: {\cal T}^h(S) \times {\Bbb Q}{\cal L}^d(S) \rightarrow
{\Bbb R}$, which is called {\em a length of a lamination} $f$ w.r.t. the complex
structure $m$.

Analogous function can be defined if we have a decorated surface $m \in {\cal
T}^d(S)$ and an unbounded rational measured lamination $f \in {\Bbb Q}{\cal
L}^h(S)$. The curves still can be transformed into geodesics, but in this case
they can have infinite length. Now, as while considering coordinates on ${\cal
T}^d$, take the distance between intersection points of the geodesics and the
horocycles around their endpoints (with negative sign if the horocycles
intersect) and take the weighted sum of these numbers. We have obtained a
function $l_{{\cal L}{\cal T}}: {\Bbb Q}{\cal L}^h(S) \times {\cal T}^d(S)
\rightarrow {\Bbb R}$, which we shall denote by the same letter and call by the
same name.

There is also the third function $l_{{\cal L}{\cal L}}:  {\Bbb Q}{\cal L}^d(S)
\times {\Bbb Q}{\cal L}^h(S) \rightarrow {\Bbb Q}$ which is called {\em an
intersection index} and is defined as follows.  Take two laminations from
${\Bbb Q}{\cal L}^d(S)$ and ${\Bbb Q}{\cal L}^h(S)$, respectively, and draw them
on $S$ in such a way, that the number of intersection points is as low as
possible.  Then the intersection index is the sum over all intersection points
of product of weighs of the intersecting curves.

The main properties of these functions are

\vspace{3mm}{\bf 1}(continuity). The functions $l_{{\cal T}{\cal L}},l_{{\cal L}{\cal T}}$
and
$l_{{\cal L}{\cal L}}$ are continuous.

\vspace{3mm}{\bf 2}(homogeneity).
\begin{equation}
l_{{\cal T}{\cal L}}(m,C f) = C l_{{\cal T}{\cal L}}(m,f),
\end{equation}
\begin{equation}
l_{{\cal L}{\cal T}}(C f,m) = C l_{{\cal L}{\cal T}}(f,m),
\end{equation}
\begin{equation}
l_{{\cal L}{\cal L}}(C f_1,f_2) = l_{{\cal L}{\cal L}}(f_1,C f_2) =
C l_{{\cal L}{\cal L}}(f_1,f_2)
\end{equation}
for any nonnegative real number $C$.

\vspace{3mm}{\bf 3}(asymptotic compatibility).
\begin{equation}
\lim_{C \rightarrow \infty} l_{{\cal T}{\cal L}}(Cm_1,m_2)/C =
l_{{\cal L}{\cal L}}(m_1,m_2) =
\lim_{C \rightarrow \infty} l_{{\cal L}{\cal T}}(m_1,Cm_2)/C.
\end{equation}
Here we have identified the spaces ${\cal L}^h(S)$ with ${\cal T}^h(S)$
and ${\cal L}^d(S)$ with ${\cal T}^d(S)$ using a graph coordinate
system. (In particular this identification gives sense to the multiplication
of a complex structure by a number.) The statement means that the equality is
true for any coordinate system.

\vspace{3mm}{\bf 4}(compatibility with coverings).
Let $\sigma: \tilde{S} \rightarrow S$ is an unramified $N$-fold covering then
\begin{eqnarray}
l_{{\cal T}{\cal T}}(\sigma^* m_1,\sigma^* m_2) &=&
N l_{{\cal T}{\cal T}}(m_1,m_2) \\
l_{{\cal L}{\cal T}}(\sigma^* f,\sigma^* m) &=&
N l_{{\cal L}{\cal T}}(f,m) \\
l_{{\cal L}{\cal L}}(\sigma^* f_1,\sigma^* f_2) &=&
N l_{{\cal L}{\cal L}}(f_1,f_2) \\
\end{eqnarray}

{\em Proof of the continuity.} Let us first prove the continuity of the
function
$l_{{\cal T}{\cal L}}$. To do this it suffices to prove it for laminations
without curves with negative weights. Indeed, if we add such curve to a
lamination the length is obviously changes continuously. Now we are going to
show that the length of integral lamination is a convex function of its
coordinates, i.e. that

\begin{equation}
l_{{\cal T}{\cal L}}(m,f_1) + l_{{\cal T}{\cal L}}(m,f_2)
\leq l_{{\cal T}{\cal L}}(m,f_1 +_\Gamma f_2), \label{ineq}
\end{equation}
were by $f_1 +_\Gamma f_2$ we mean a lamination with coordinates being sums of the
respective coordinates of $f_1$ and $f_2$. Taking into account the homogeneity
property of $l_{{\cal T}{\cal L}}$, one sees that the inequality (\ref{ineq})
holds for all rational laminations and therefore can be extended by continuity
for all real laminations.

Prove now the inequality (\ref{ineq}).  Draw both laminations $f_1$ and $f_2$
on the surface and deform them to be geodesic.  These laminations in general
intersect each other in finite number of points.  Then retract the whole
picture to the fat graph in such a way that no more intersection points appear and the
existing ones are moved to the edges.  Now it becomes obvious, that at each
intersection point we can rearrange our lamination cutting both intersecting
curves at the intersection point and gluing them back in another order in a way
to make the resulting set of curves homotopically equivalent to a
nonintersecting collection.  We can do it at each intersection point in two
different ways, and we use the retraction to the fat graph to choose one of
them.

Indeed, connect them as shown on fig. \ref{convex}.

One can easily see that the
numbers on edges, corresponding to the new lamination $f$ are exactly the sums
of the numbers, corresponding to $f_1$ and $f_2$. On the other hand, the
lamination on the original surface is no longer geodesic, because the curves
of it have breaks. But its length is exactly $l_{{\cal T}{\cal L}}(f_1,m) +
l_{{\cal T}{\cal L}}(f_2,m)$. When we deform $f$ to a geodesic lamination,
its length can only decrease, what proves the inequality (\ref{ineq}).

\vspace{5mm}
\epsfxsize9cm
\centerline{\epsfbox{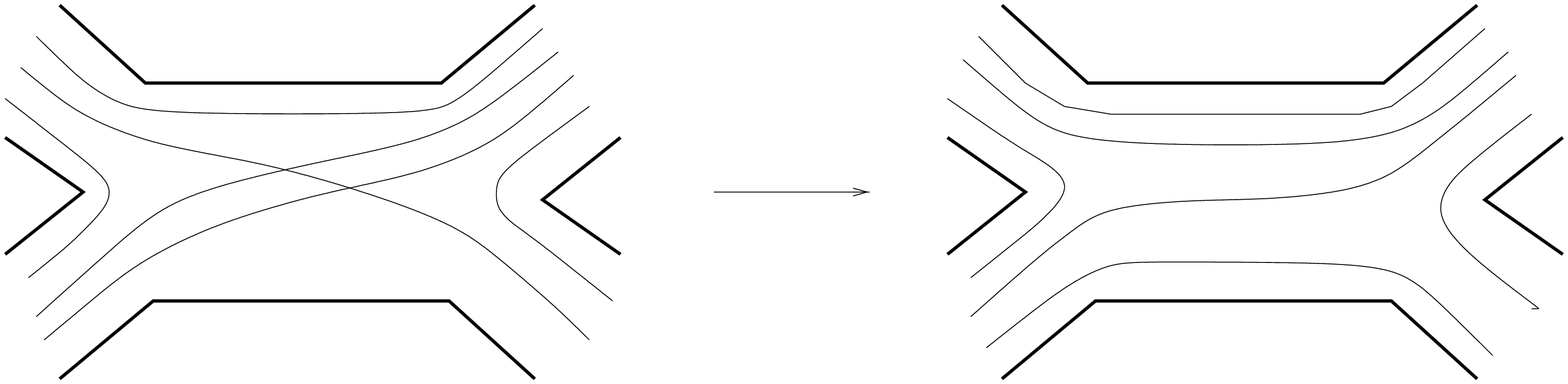}} 
\vspace{5mm}
\refstepcounter{equation}\label{convex}
\vspace{5mm}
\centerline{Fig. \arabic{equation}}

Proves of the continuity of $l_{{\cal L}{\cal T}}$ and $l_{{\cal L}{\cal L}}$
are absolutely analogous and we do not repeat them here. Note only that the
proof of the continuity of the length of unbounded laminations is even
simpler since we don't need to consider the case with curves with negative
weights separately.

Note that given a graph $\Gamma \in \Gamma_0^3(S)$  the functions $l_{{\cal
L}{\cal L}}$ as well as $l_{{\cal L}{\cal T}}$ are given by very simple
formulae provided the coordinates of the unbounded lamination are
nonnegative.
\begin{equation} \label{LL}
l_{{\cal L}{\cal L}}(f^H(z_1,\ldots,z_N),f^d(v_1,\ldots,v_N))  =
\sum_{\alpha \in E(\Gamma)} z_\alpha v_\alpha;
\end{equation}
\begin{equation} \label{LT}
l_{{\cal L}{\cal T}}(f^H(z_1,\ldots,z_N),m^d(u_1,\ldots,u_N)) =
\sum_{\alpha \in E(\Gamma)} z_\alpha u_\alpha,
\end{equation}
where $f^H(z_1,\ldots,z_N),f^d(v_1,\ldots,v_N)$ and $m^d(u_1,\ldots,u_N)$ are an unbounded
lamination, a bounded lamination and a decorated surface given by the
respective coordinates w.r.t. $\Gamma$. Taking into account independence on the
choice of $\Gamma$ and continuity of the functions $l_{{\cal L}{\cal L}}$ and
$l_{{\cal L}{\cal T}}$ and the fact that for almost all unbounded laminations
one can make
all coordinates on edges
positive by changing $\Gamma$ and orientations of the holes, we can in
principle compute these functions for any values of the arguments.

The coincidence of the r.h.s.  of these formulae immediately gives the proof of
the second part or the asymptotic compatibility property.  To demonstrate the
first part it is sufficient to check it for the bounded lamination being a
closed curve $\gamma$ of weight one.  Without loss of generality we may assume
that we have a coordinate system where the coordinates of the holed surface are
positive.  Instead of computing the length $l$  of the geodesics, compute the
function $2\cosh (l/2) = {\rm Tr}(M(\gamma))$, where $M(\gamma)$ is
the element of the Fuchsian group $\Delta$ corresponding to $\gamma$.  This
function has the same leading term in $C$.  One computes the trace using the
construction of the Fuchsian group as a trace of product of matrices
$A(z_\alpha)$ and $I$ from the section 5.1.  But if all values of coordinates
are big and positive, we can replace the matrices $A(z_\alpha)$ by
$e^{z_\alpha/2} \left(\begin{array}{cc} 0 & 1\\ 0 & 0 \end{array} \right)$
without changing the leading term.  It is also easy to see that after such
replacement the trace of the product of matrices along a curve is always
nonzero (here we need the positivity of the coordinates) and proportional to
the exponent of the r.h.s.  of (\ref{LL}).  Note that we have used nowhere in
the proof that the curve has no self intersections.

\section{Weil-Petersson forms.}
  Since the holed Teichm\" uller space ${\cal T}^h(S)$ is a subspace of the
space of representations of the fundamental group of the surface in
$PSL(2,{\Bbb R})$, it possesses a canonical Poisson structure called {\em the
Weil-Petersson structure} (one of possible constructions of the Poisson structures
on representation spaces of fundamental groups of two dimensional surfaces with
proves of some properties is presented in \cite{FR}). This
Poisson structure is degenerate and its Casimir functions are the lengths of
geodesics surrounding holes. Therefore, if we fix all such lengths we get a
symplectic leaf. In particular the space ${\cal T}(S)$ is a symplectic leaf
as a subspace of ${\cal T}^h(S)$ of the surfaces with vanishing lengths
of geodesics surrounding holes.

On the other hand the space ${\cal T}^d(S)$ projects onto ${\cal T}(S)$.  So
we can invert the Poisson structure on ${\cal T}(S)$ and pull back the
resulting symplectic structure to ${\cal T}^d(S)$.  The resulting degenerate
two-form is also called the {\em Weil-Petersson form}.

It turns out that these forms have very simple expressions in terms of the
constructed coordinates:
\begin{equation} \label{owp}
\omega_{WP} = \sum_{\alpha \in EE(\Gamma)}du_\alpha \wedge du_{\alpha(1)}
\end{equation}
See \cite{Penner1} for the proof.
\begin{equation}  \label{pwp}
P_{WP} = \sum_{\alpha \in EE(\Gamma)}\frac{\partial}{\partial z_\alpha} \wedge
\frac{\partial}{\partial z_{\alpha(1)}}
\end{equation}
The proof is given in \cite{FZ}.

One can observe (which was done in \cite{Papadop} for decorated surfaces)
that these formulae give the Poisson structure and degenerate 2-form on the
respective spaces of laminations. Although the proofs that the formulae give
the Weil-Petersson forms indeed requires relatively long computations one can
easily check the following properties of the forms $\omega_{WP}$ and $P_{WP}$.

\vspace{3mm}{\bf 1.} The expressions for $\omega_{WP}$ and $P_{WP}$ give
the same forms independently on the graph $\Gamma$.

{\bf 2.} If $\sigma:\tilde{S} \rightarrow S$ is an unramified $N$-fold
covering let $h_1,h_2:{\cal T}^H(\tilde{S})\rightarrow {\Bbb R}$ be
two functions generating Hamiltonian vector field tangent to
$\sigma^*{\cal T}^H(S)$. (The latter condition is equivalent to the demand
that for any function $h:{\cal T}^H(\tilde{S})\rightarrow {\Bbb R}$ such
that $\sigma^{**}h = 0$ we have $\{h_1,h\} = \{h_2,h\} = 0$.)
Then
\begin{equation}
\sigma^{**}\{h_1,h_2\} = N\{\sigma^{**}h_1,\sigma^{**}h_1\}.
\end{equation}
$\{,\}$ denote the Poisson bracket given by $P_{WP}$.

As a sketch of a proof of the formula (\ref{pwp}) note that it is enough
to check that it gives the correct (e.g. taken from Goldman's paper
\cite{Goldman} ) Poisson bracket for, say, lengths of any two closed
geodesics intersecting at one point. Then taking into account the property
{\bf 2} one can extend this equality for Poisson brackets between lengths of
any two geodesics and therefore prove that (\ref{pwp}) gives indeed the
Weil-Petersson Poisson bracket.

\section{Applications and remarks.}
In this section we are describing different subjects related to our
description of Teichm\" uller and lamination spaces. The work on these
subjects is still in progress and we hope to publish more detailed text on
them in the nearest future. Nevertheless we decided to write down some clear
parts of these subjects in order to convince the reader that the considered
approach has several amazing perspectives.

\subsection{Action of the mapping class group.}

Here we shall briefly discuss the properties of the action of the mapping
class group ${\cal D}(S)$ on the lamination and the Teichm\" uller spaces.

  The action of ${\cal D}(S)$ on the space of bounded laminations
${\cal L}^d(S)$ is nowhere discontinuous. Indeed, a generic rational
lamination is a collection of $3g-3+s$ curves. If we forget
about weights there is only finite number of ${\cal D}(S)$-orbits of such
collections. Therefore any ${\cal D}(S)$-orbit of a rational lamination
intersects with a submanifold of dimension $3g-3+s$. Since
${\rm dim}{\cal L}^d(S)= 6g-6+3s$ the quotient of any open dense
subset of ${\cal L}^d(S)$ by ${\cal D}(S)$ is nonhausdorf and even a
space without closed points other than zero.

  The quotient of ${\cal L}^H(S)$ by ${\cal D}(S)$ is nonhausdorf as well
however there exists an open dense subset of it with a Hausdorf quotient.
Indeed let ${\cal L}^H_0(S) \subset {\cal L}^H(S)$ be the set of laminations
supported on $6g-6+3s$ nonclosed curves. This space is dense, since in any
neighbourhood of any rational lamination there obviously exists a lamination
$f \in {\cal L}^H_0(S)$ and open, since we span a neighbourhood of $f$ changing
weights on these curves. On the other hand the curves of such laminations
cut the surface into cells and a graph dual to this decomposition gives us
the canonical graph $\Gamma(f) \in \Gamma^3_0(S)$. If orientations of all
holes are induced from the orientation of the surface all coordinates of
$f$ w.r.t. $\Gamma$ are positive. Inversely, if all coordinates of a
lamination $f$ w.r.t. a graph $\Gamma \in \Gamma^3_0(S)$ are positive
then $f \in {\cal L}^H_0(S)$ and $\Gamma$ is just the graph associated
to $f$.

  Equivalently we can say that almost any lamination $f \in {\cal L}^H_0(S)$
given by graph coordinates can be transformed by changing the graph and
orientations of holes to a lamination with positive coordinates and the
final graph and coordinates are uniquely defined by $f$ up to the final
graph symmetry.

  The actions of the mapping class group on ${\cal T}^H(S)$ and
${\cal T}^d(S)$ are properly discontinuous and the quotients are
well defined Hausdorf spaces and even orbifolds ${\cal M}^H(S)$
and ${\cal M}^d(S)$. The fundamental domains of the ${\cal D}(S)$-action
on ${\cal T}^d(S)$ are described by R. Penner in \cite{Penner} just in terms
of graph coordinates. The fundamental domains of the ${\cal D}(S)$-action
on ${\cal T}^h(S)$ are described by S.Kojima \cite{Kojima}. However, although
it is possible to describe his domains in graph coordinates, we do not know
any simple expression for them. Note, that the constructions of Penner and
Kojima give not only the fundamental domains but also a full cell
decomposition of the spaces ${\cal M}^H(S)$ and ${\cal M}^d(S)$,
respectively.

\subsection{Quantisation.}

  Once the graph $\Gamma$ (or in another words the
decomposition of the surface $S$ into triangles) is chosen on can easily
quantise the corresponding Teichm\" uller space ${\cal T}^H(S)$ in the
following sense.  Consider the $*$-algebra ${\cal T}^\hbar(\Gamma)$ generated
by real generators $\{Z^\hbar_\alpha | \alpha \in E(\Gamma)$ ({\em real} mean
that $(Z^\hbar_\alpha)^* = Z^\hbar_\alpha$) with relations
\begin{equation} \label{comm}
[Z^\hbar_\alpha, Z^\hbar_\beta ] = 2\pi i\hbar\{Z_\alpha, Z_\beta\}
\end{equation}

 This algebra has an obvious center generated by
\begin{equation} \label{center}
\{P_\gamma| \gamma \in F(\Gamma),\
P_{\gamma} = \sum_{\alpha \in \gamma} Z^\hbar_\alpha\}.
\end{equation}

  It is not a big deal to  describe all irreducible $*$-representations
of this algebra using Stone--von Neumann theorem. An irreducible representation
is unambiguously fixed by the values of the operators $P_\gamma$ for all
$\gamma in F(\Gamma)$ which must be scalars.  For example one can represent
all operators $Z^\hbar_\alpha$ in $L^2({\Bbb R}^n)$, where
$n = \frac{1}{2}( \sharp E(\gamma) - \sharp F(\Gamma))$, by linear
combinations with real coefficients of constants and the operators $x_i$
and $i\frac{\partial}{\partial x_i}$, where
$\{x_i|i=1,\ldots,n\}$ is a standard coordinate system on ${\Bbb R}^n$.
Now our task is to identify the $*$-algebras constructed using different
graphs $\Gamma$ and $\Gamma^\prime$ corresponding to a given surface $S$ .
In order to make this identification we just construct a $*$-homomorphism
$K(\Gamma, \Gamma^\prime): {\cal T}^\hbar(\Gamma) \rightarrow {\cal
T}^\hbar(\Gamma^\prime)$ of the $*$-algebra generated by
$\{Z_{\alpha^\prime}|\alpha^\prime \in E(\Gamma^ \prime)\}$ to the algebra
generated by  $\{Z^\hbar_\alpha|\alpha \in E(\Gamma )\}$.  We require this
homomorphism to have the following properties:

1. {\em Classical limit}. We demand that the algebra homomorphism should
tend to the classical homomorphism of the algebras fo function on ${\cal
T}^H$ when the parameter $\hbar$ tends to zero.

2. {\em Path independence}. We demand that if we have three graphs
$\Gamma$, $\Gamma^\prime$ and $\Gamma^{\prime \prime}$ then the
homomorphisms should satisfy the condition
$K(\Gamma^{\prime \prime}, \Gamma^\prime)K(\Gamma^\prime, \Gamma) =
K(\Gamma^{\prime \prime}, \Gamma)$.

  Using the latter demand one can obviously reconstruct the homomorphism
$K(\Gamma^\prime, \Gamma)$ for any $\Gamma$ and $\Gamma^\prime$ once one
knows these homomorphisms for pairs of graphs related by single flips.
However an arbitrary set of such flip homomorphisms {\it a priori} does
not satisfy  the path independence condition since one can get one graph from
another by different sequences of flips.  Or in another words one must check
that if a sequence of flips does not change the graph then the corresponding
product of algebra homomorphisms is the identical one.

Fortunately one can check the latter condition  only for one sequence of
flips since others  are just compositions of this one.  (This is a kind of
folklore statement.  I would be grateful to anybody who
let me know a reference with a nice short proof of it) \footnote{I am
indebted to R. Lawrence who pointed me out
that this statement belongs to McLane}.

\noindent
\begin{minipage}{9cm}
~~~Describe now this distinguished sequence of flips. Consider  two edges
having exactly one common vertex. One can easily see that a sequence of
five flips of these edges (such that we never flip the same edge twice
consequtively) does not change  the graph. It is may be more geometrically
transparent to see this on the dual graph where the two edges correspond to
two edges separating three triangles forming a pentagon. A pentagon can be
cut into three triangles in only five possible ways which are related by
flips (fig.\ref{pent1})
\end{minipage}
~~~~~
\begin{minipage}{5cm}
\unitlength 0.2mm
\centerline{
\begin{picture}(0,0)(100,-145)
\put(0,0){\line(2,1){ 40}}
\put(40,20){\line(2,-1){ 40}}
\put(80,0){\line(-1,-2){ 20}}
\put(60,-40){\line(-1,0){40}}
\put(20,-40){\line(-1,2){ 20}}
\put(0,0){\line(1,0){80}}
\put(0,0){\line(3,-2){60}}
\end{picture}
\begin{picture}(0,0)(20,-200)
\put(0,0){\line(2,1){ 40}}
\put(40,20){\line(2,-1){ 40}}
\put(80,0){\line(-1,-2){ 20}}
\put(60,-40){\line(-1,0){40}}
\put(20,-40){\line(-1,2){ 20}}
\put(0,0){\line(3,-2){60}}
\put(40,20){\line(1,-3){20}}
\end{picture}
\begin{picture}(0,0)(-60,-145)
\put(0,0){\line(2,1){ 40}}
\put(40,20){\line(2,-1){ 40}}
\put(80,0){\line(-1,-2){ 20}}
\put(60,-40){\line(-1,0){40}}
\put(20,-40){\line(-1,2){ 20}}
\put(40,20){\line(1,-3){20}}
\put(40,20){\line(-1,-3){20}}
\end{picture}
\begin{picture}(0,0)(-30,-65)
\put(0,0){\line(2,1){ 40}}
\put(40,20){\line(2,-1){ 40}}
\put(80,0){\line(-1,-2){ 20}}
\put(60,-40){\line(-1,0){40}}
\put(20,-40){\line(-1,2){ 20}}
\put(40,20){\line(-1,-3){20}}
\put(20,-40){\line(3,2){60}}
\end{picture}
\begin{picture}(0,230)(70,-65)
\put(0,0){\line(2,1){ 40}}
\put(40,20){\line(2,-1){ 40}}
\put(80,0){\line(-1,-2){ 20}}
\put(60,-40){\line(-1,0){40}}
\put(20,-40){\line(-1,2){ 20}}
\put(20,-40){\line(3,2){60}}
\put(0,0){\line(1,0){ 80}}
\end{picture}
}

\refstepcounter{equation}\label{pent1}
\vspace{5mm}
\centerline{Fig. \arabic{equation}}
\end{minipage} 

Now we are going to give an answer for the flip homomorphism satisfying the
above conditions.
The simplest way to describe this rule is to draw the picture (\ref{rule}).
On the left picture a fragment of the graph $\Gamma$ is shown together with
the algebra elements associated to the edges. On the right hand picture the
corresponding fragment of the graph $\Gamma^\prime$ together with  the
operators images of operators corresponding to the edges in the former
algebra. The remainder of the graphs as well as the operators assigned
to the remainder of the edges of $\Gamma$ and $\Gamma^\prime$ coincide.

\centerline{
\setlength{\unitlength}{1.8mm}%
\begin{picture}(50,29)(2,48)
\thicklines
\put(28,70){\line( 1,-2){ 4}}
\put(32,62){\line( 1, 0){28}}
\put(60,62){\line( 1, 2){ 4}}
\put(60,62){\line( 1,-2){ 4}}
\put(32,62){\line(-1,-2){ 4}}
\thinlines
\put(18,62){\vector(-1, 0){  0}}
\put(18,62){\vector( 1, 0){ 5}}
\thicklines
\put(10,54){\line( 2,-1){ 8}}
\put(10,70){\line( 0,-1){16}}
\put(10,54){\line(-2,-1){ 8}}
\put( 2,74){\line( 2,-1){ 8}}
\put(10,70){\line( 2, 1){ 8}}
\put( 4,74){\makebox(0,0)[lb]{$A^\hbar$}}
\put(16,74){\makebox(0,0)[rb]{$B^\hbar$}}
\put(12,62){\makebox(0,0)[lc]{$Z^\hbar$}}
\put(16,50){\makebox(0,0)[rt]{$C^\hbar$}}
\put( 4,50){\makebox(0,0)[lt]{$D^\hbar$}}
\put(30,54){\makebox(0,0)[lt]{$D^\hbar - \phi(-Z^\hbar)$}}
\put(62,54){\makebox(0,0)[rt]{$C^\hbar+\phi(Z^\hbar)$}}
\put(62,69){\makebox(0,0)[rb]{$B^\hbar-\phi(-Z\hbar)$}}
\put(30,69){\makebox(0,0)[lb]{$A^\hbar+\phi(Z^\hbar)$}}
\put(47,64){\makebox(0,0)[cb]{$-Z^\hbar$}}
\end{picture}
}
\refstepcounter{equation}\label{rule}
\vspace{5mm}
\centerline{Fig. \arabic{equation}}
\vspace{5mm}

where $\phi$ is a real-analytic function of
$x$ depending on $\hbar$ as on a parameter
\begin{equation} \label{phi}
\phi(x) =
-\frac{\pi\hbar}{2}\int_{\Omega} \frac{e^{-ipx}}{\sinh(\pi p)\sinh(\pi \hbar
p)}dp,
\end{equation}
and the contour $\Omega$ is the real axis shifted slightly to the
upper half plane at the origin.

The constructed isomorphisms show that in fact the algebra ${\cal
T}^\hbar(\Gamma)$ does not depend on a particular choice of the graph
$\Gamma$ and we redenote this algebra by ${\cal T}^\hbar(S)$

The construction of quantisation alowes to make the following
constructions, statements, conjectures and remarks.

\vspace{3mm}{\bf 1.} {\em Projective unitary mapping class group
representations.} (construction and a statement)

  The above construction gives us representations of the mapping
class group ${\cal D}(S)$ of a nonclosed surface $S$ in a suitably completed
Heisenberg $*$-algebra. Our aim now is to construct a unitary projective
representations of certain subgroups of ${\cal D}(S)$ in a Hilbert space.

Assign real numbers $l_1,\ldots,l_s$ to the holes of the
surface $S$. Let now $H(S,l_1,\ldots,l_s)$ be the Hilbert space of the
unitary representation of the Heisenberg algebra ${\cal T}^\hbar(S)$, where
$l_1,\ldots,l_s$  are the values of the central elements (\ref{center})
corresponding to the holes. The mapping class group ${\cal D}(S)$ obviously
permute the numbers $l_1,\ldots,l_s$. Let ${\cal D}(S,l_1,\ldots,l_s) \subset
{\cal D}(S)$ be the stabiliser of the collection $l_1,\ldots,l_s$.

The following proposition follows obviously from the classification of the
algebra ${\cal T}^\hbar(S)$ representations.

For a given value of the constant $\hbar$ There exists a central extension
$\tilde{\cal D}^\hbar(S,l_1,\ldots,l_s)$ of the group ${\cal
D}(S,l_1,\ldots,l_s)$ and a unique unitary representation
$T^\hbar(S,l_1,\ldots,l_s)$ of this group in $H(S,l_1,\ldots,l_s)$ such that
for any elements $x \in \tilde{{\cal D}(S,l_1,\ldots,l_s)}$ and $a \in {\cal
H}^\hbar(S)$
\begin{equation}
T(a)x = a(x)T(a),
\end{equation}
where both
sides are elements of ${\mbox{End}}H(S,l_1,\ldots,l_s)$ and $a(x)$ is the
result of the action of the mapping class group element $a$ on $x$
constructed above.

Since we are not gotng to discuss the central extension $\tilde{\cal D}$ in
details we shall call the above described representations as projective
representations of ${\cal D}$ instead of ordinary ones of ${\cal D}$ just in
order to simplify notations.

\vspace{3mm}{\bf 2.}{\em $\left(\hbar \rightarrow
\frac{1}{\hbar}\right)$-invariance.} (statement)

  The algebras generated by $\{Z^\hbar_\alpha\}$ corresponding to different
values of $\hbar \geq 0$ are obviously isomorphic to each other. However these
isomorphisms are in general not canonical. It means that {\em a priory} it is
not possible to define an isomorphism commuting with the action of the
mapping class group.

It turns out that for some pairs of values of the parameter $\hbar$ such
equivariant isomorphism does exist. Indeed consider the $*$-algebra
isomorphism $D(\Gamma): {\cal T}^\hbar(\Gamma) \rightarrow
{\cal T}^\frac{1}{\hbar}(\Gamma)$ given on generators by
\begin{equation} \label{dual}
D(\Gamma): Z^\hbar_\alpha \mapsto \frac{1}{\hbar}Z^\frac{1}{\hbar}_\alpha
\end{equation}
One can easily check that for two different graphs $\Gamma$ and
$\Gamma^\prime$ related by a flip one has
\begin{equation}
K(\Gamma, \Gamma^\prime) D(\Gamma) = D(\Gamma^\prime) K(\Gamma,
\Gamma^\prime).   \label{commut}
\end{equation}
what means that in fact this isomorphism does not depend on the graph we have
chosen to define it. For this reason we shall denote it below as $D(S)$.

The proof of the equality (\ref{commut}) is straightforward provided one uses
the property verified by the function $\phi(x,\hbar)$:
\begin{equation}
\phi(x,\hbar) = \hbar\phi(\frac{1}{\hbar}x,\frac{1}{\hbar}).
\end{equation}

This property shows in particular the isomorphism between representations of
the mapping class group
\begin{equation}
T^\frac{1}{\hbar}(S,\frac{l_1}{\hbar},\ldots,\frac{l_s}{\hbar}) \cong
T^\hbar(S,l_1,\ldots,l_s).
\end{equation}

It means also that the
quantisation has two isomorphic classical limits $\hbar \rightarrow 0$ and
$\hbar \rightarrow \infty$. If we assume that we have here the quantum
Liouville theory we can interpret the the parameter $\hbar$ as a coupling
constant and the $\hbar \leftrightarrow \frac{1}{\hbar}$ symmetry as a
"week--strong coupling constant duality".

\vspace{3mm}{\bf 3.}{\em Many more unitary projective mapping class group
representations.} (a construction and a statement)

 The family of representations constructed in the previous
remark give for free a wide class of representations (also projective and
unitary).  Indeed, consider an unramified $N$-fold covering $\sigma:
\tilde{S} \rightarrow S$.  Consider the representation
$T^\hbar(\tilde{S},\tilde{l}_1,\ldots, \tilde{l}_{\tilde{n}})$ of the mapping
class group ${\cal D}(\tilde{S})$.  Restrict this representation to the
congruence subgroup ${\cal D}(S,\sigma)$ and then induce the restriction to
the whole group ${\cal D}$. Denote the resulting representation by
$T^\hbar(S,\sigma,\tilde{l}_1,\ldots, \tilde{l}_{\tilde{n}})$ although it
depends only on the ${\cal D}(S)$-orbit of $\sigma$ in the space of coverings
of $S$.

The decomposition of the constructed representations is yet unclear for us.
however one can check the following property:

  There exists a canonical mapping
\begin{equation}
T^\frac{\hbar}{N}(S,l_1,\ldots, l_n)
T^\hbar(S,\sigma,\tilde{l}_1,\ldots, \tilde{l}_{\tilde{n}})
\end{equation}

\vspace{3mm}{\bf 4.}{\em Geodesic length operators.}(a construction and a statement)
  For any closed unoriented path $\gamma$ on $S$ one can associate a smooth
real function $l_\gamma$ on ${\cal T}^H(S)$. $l_\gamma$ is the
length of a closed geodesic in the homotopy class of $\gamma$. Introduce also
another set of functions $L_\gamma$ for each path $\gamma$ just as $L_\gamma
= 2\cosh l_\gamma$.

  The functions $\{L_\gamma\}$ were studied by several authors and mainly by
W.M.Goldman \cite{Goldman}. Among their nice properties let us mention here
that they generate a Poisson algebra (w.r.t. the multiplication and the
Weil-Petersson Poisson bracket) over ${\Bbb Z}$. (It means that a product and
a Poisson bracket of two such functions is a linear combination of such
functions with integral coefficients.) The quantum deformation of this
algebra is also known \cite{Turaev}.

  The aim of this remark is to embed this algebra into a suitable
completion of the constructed algebra ${\cal T}^\hbar(S)$.

The function $\L_\gamma$ can be easily expressed for any $\gamma$ in terms of
graph coordinates on ${\cal T}^H$. For any $\gamma$ it is given by an
expression of the form:
\begin{equation}
L_\gamma = \sum_{j\in J}e^{\frac{1}{2} \sum_{\alpha \in E(\Gamma)}
m_j(\gamma,\alpha) z_\alpha},  \label{clen}
\end{equation}
where $m_j(\gamma,\alpha)$ are
certain integral numbers and $J$ is just  a finite set of indices
numerating the terms in (\ref{clen}).

Let us now define and formulate some properties of quantum analogues of these
functions.

Denote by $\widehat{\cal T}^\hbar$ a completion of the algebra ${\cal
T}^\hbar$ containing $e^{xZ_\alpha}$ for any real $x$.

Let for any closed path $\gamma$ on $S$ such that it never goes along the same edge twice the operator $L^\hbar \in
\hat{\cal T}^\hbar$ is given by
\begin{equation} \label{qlen}
L^\hbar_\gamma = \sum_{j\in J}e^{\frac{1}{2} \sum_{\alpha \in E(\Gamma)}
m_j(\gamma,\alpha) Z^\hbar_\alpha},
\end{equation}
where the numbers $m_j(\gamma,\alpha)$ are the same as in (\ref{clen}).

Note that the operators $\{L^\frac{1}{\hbar}_\gamma\}$ can be considered as
belonging to the algebra $\hat{\cal T}^\hbar$ due to the isomorphism
$D(S)$ given by (\ref{dual}). In terms of the generators of $\hat{\cal
T}^\hbar$ they are obviously given by
\begin{equation} \label{dlen}
L^\frac{1}{\hbar}_\gamma = \sum_{j\in J}e^{\frac{1}{2\hbar} \sum_{\alpha \in
E(\Gamma)} m_j(\gamma,\alpha) Z^\hbar_\alpha},
\end{equation}.
Unfortunately this definition is not good for curves going along an edge two or more times. However is is sufficient to define operators for all curves since for any curve one can make it go along each edge no more then once by changing the graph.

The properties of the operators $\{L^\hbar_\gamma\}$ and
$\{L^\frac{1}{\hbar}_\gamma\}$ we would like to mention here are the following:

{\em 1).} The operators $\{L^\hbar_\gamma\}$ and
$\{L^\frac{1}{\hbar}_\gamma\}$ are correctly defined, what means that they depend only on the homotopy class of the path $\gamma$.

{\em 2).} For any $\gamma$ and $\gamma^\prime$ the operators $L^\hbar_\gamma$
and  $L^\frac{1}{\hbar}_{\gamma^\prime}$ commute.

{\em 3).} If two closed paths $\gamma$ and $\gamma^\prime$ do not intersect
then the operators $L^\hbar_\gamma$ and  $L^\hbar_{\gamma^\prime}$
commute.

{\em 4).} The algebra generated by $L^\hbar_\gamma$ (resp., $L^\frac{1}{\hbar}_\gamma$ is isomopphic to the Turaev quantum loop algebra \cite{Turaev} for the deformation parameter $q=e^{2\pi i \hbar}$ (resp., $\tilde{q} = e^{2\pi i \frac{1}{\hbar}}$). Tn particular a product $L^\hbar_\gamma L^\hbar_{\gamma^\prime}$ for any two paths $\gamma$ and $\gamma^\prime$ is a linear combination of functions $L^\hbar_{\gamma_i}$, where $\{\gamma_i\}$ is a finite set of curves, the coefficients being Laurent polynomials in $q$ with positive integer coefficients. The same is true of course for the  algebra generated by $L^\frac{1}{\hbar}_\gamma$.

After formulating such nice properties of the operators $L^\hbar_\gamma$ we
should mention some of their properties which strongly reduce their
applicability to our problems. We are studying the $*$-representations of the
algebra ${\cal T}^\hbar$ where $Z^\hbar_\alpha$ are represented by unbounded
self-adjoint operators. The exponents of such operators  are not good
operators in a Hilbert space. In particular they are defined only on
functions which can be analytically extended to a certain domain around the
real axis and/or have certain exponential decrease at infinity. However they
are still useful as a tool to study the mapping class group representations
(in particular in finding Dehn twists spectra) since the difference equations
are much simpler than the integral ones.

{\bf 5.}  {\em Modular functor} (conjecture). The association
$(S,l_1,\ldots,l_n) \rightarrow H(S,l_1,\ldots,l_n)$ is a unitary modular
functor. It means that the constructed mapping class group representations is
compatible with embeddings of one surface into another in the following way.

  Let $S^1 \rightarrow S^2$ is an embedding of surfaces. Then obviously there
exists a canonical embedding of mapping class groups ${\cal D}(S^1)
\rightarrow {\cal D}(S^2)$. Then the modular functor condition means that if
we restrict the representation of $T^\hbar(S_2,l^2_1,\ldots,l^2_{n_1})$ to
the subgroup ${\cal D}(S^1)$ we get in the decomposition to irreducible
representations the representations $T^\hbar(S_2,l^1_1,\ldots,l^1_{n_1})$
for different values of $l^1_1,\ldots,l^1_{n_1}$ only. In particular if $S^1$
is obtained form $S^2$ by cutting along a simple closed curve $\gamma$ then
\begin{equation}
T^\hbar(S_2,l^2_1,\ldots,l^2_{n_1}) =
\int_{-\infty}^{\infty} T^\hbar(S^1,l^2_1,\ldots,l^2_{n_1},l,-l)dl.
\end{equation}

{\bf 6.}  {\em Irreducibility} (conjecture). The representations
$T^\hbar(S,l_1,\ldots,l_{n_1})$ of the mapping class group are irreducible
for irrational values of the parameter $\hbar$.

\subsection{Markov numbers.}

 Consider a torus with one hole $T$.  The space of homotopy classes of simple
(i.e., without intersections) unoriented closed paths on it can be
parameterised by points of ${\Bbb Q}P^1$. Indeed, once we have chosen an
orientation of the path, we can consider it as an element of the first
homology of $T$ with compact support. It is also obvious that any simple
(indivisible) class is represented by a unique simple oriented closed path.
Since the first homology group is ${\Bbb Z}^2$, it just gives the desired
parameterisation.

Introduce the equiharmonic complex structure on $T$, i.e. the structure
which has maximal symmetry group ${\Bbb Z}/3{\Bbb Z}$. For any closed
path $\gamma$ on $T$ without self-intersections the numbers $X_\gamma =
\frac{2}{3}\cosh l(\gamma)$, where $l(\gamma)$ are the lengths of the
corresponding geodesics, are called {\em Markov numbers}.

The main properties of the Markov numbers are the following:

\vspace{3mm}{\bf 1.} Markov numbers are positive integral.

\vspace{3mm}{\bf 2.} Markov numbers include Fibonacci numbers with even numbers
$2,5,13,34,\ldots$.

Call {\em Markov triple} a triple of Markov numbers $(X,Y,Z)$ corresponding
to three geodesics having pairwise one intersection point.

\vspace{3mm}{\bf 3.} Elements of a Markov triple satisfy the Markov equation:
\begin{equation}
X^2+Y^2+Z^2 = 3XYZ
\end{equation}

\vspace{3mm}{\bf 4.} Any integer solution of this equation is a Markov triple.

\vspace{3mm}{\bf 5.} For any Markov triple $(X,Y,Z)$ the triples $(Y,Z,X)$ and $(Z,Y-3XZ,X)$ are also Markov triples. Any Markov triple can be obtained from the triple $(1,1,1)$ by a sequence of such transformations.

Since homotopy classes of closed nonselfintersecting curves can be
parameterised by ${\Bbb Q}P^1$, one can choose an affine coordinate on ${\Bbb
Q}P^1$ in such a way that the curves with coordinates $0,1$ and $\infty$ have
Markov numbers $1$. Denote by $M(u)$ the Markov number corresponding to the
curve with the coordinate $u \in {\Bbb Q}$.

\vspace{3mm}{\bf 6.} The function $\psi(\frac{p}{q}) = \frac{1}{q}{\rm
arcosh}(\frac{3}{2}M(\frac{p}{q}))$, where $\gcd(p,q) = 1$, is extendible to
a continuous convex function on ${\Bbb R}$.

\vspace{3mm}{\bf 7.} $M(x) = M(1-x) = M(\frac{1}{x}) = M(\frac{1}{1-x}) =
M(\frac{x}{x-1}) = M(\frac{x-1}{x})$

\vspace{3mm}{\bf 8.} For any closed geodesics $\gamma$ on $S$ there exists a unique
geodesics $\gamma\prime$ going from the puncture to the puncture which do not
intersect $\gamma$. Let $l(\gamma\prime)$ be the length of the piece of
$\gamma\prime$ between the intersection points with the horocycle surrounding
area $3$. Then $e^{l(\gamma\prime)} = M(\gamma).$

\vspace{3mm}{\bf 9.} (Markov conjecture). The famous unproven Markov conjecture says that
two Markov numbers $M(x)$ and $M(y)$ are different unless $x$ and $y$ are
related by transformations from property 7.

Taking into account that the segment $[0,1]$ is the fundamental domain of the
action of transformations from property {\bf 7}, one can reformulate the
Markov conjecture as that if $M(x) = M(y)$ and $x,y \in [0,1]$ then $x=y$.

{\em Proves of the properties.} (unfortunately, without the last one and the
property 4).

There is only one graph corresponding to a holed torus. It has two vertices,
three edges and one face. This graph has obvious ${\Bbb Z}/3{\Bbb Z}$
symmetry group cyclically permuting the edges. Let $x,y,z$ be the coordinates
on the Teichm\" uller space ${\cal T}^H(S)$ w.r.t. this graph.

 A closed curve on $S$ can be considered as a bounded lamination if we assign
the weight $1$ to it. The standard graph coordinates of such laminations are
given by three nonnegative integers $n_1,n_2,n_3$. These three numbers have
no common factor, because otherwise the weight of the curve would be greater
than $1$. On the other hand one of the numbers should be a sum of two others
since otherwise there would be a component surrounding the hole. The relation
between this parameterisation by $n_1,n_2,n_3$ and the parameterisation by
${\Bbb Q}P^1$ described above is given by
\begin{equation}
x =\left\{
\begin{array}{ll} \frac{-n_2}{n_1} & {\rm if }\  n_3 = n_1+n_2\\
                      \frac{n_2}{n_1} & {\rm if }\ n_1 = n_2+n_3\   {\rm or }\
                         n_2 = n_3+n_1
\end{array} \right.
\end{equation}

Denote by ${\rm Z},{\rm X}$ and ${\rm Y}$ one thirds of traces of the
elements of the Fuchsian group corresponding to the curves with
co\-or\-di\-na\-tes $(1,1,0)$,$(1,0,1)$ and $(0,1,1)$, respectively. They can
be easily computed using the explicit formulae for the Fuchsian group:
\begin{displaymath} 
{\rm Z}=\frac{1}{3}(e^{(x+y)/2}+e^{(x-y)/2}+e^{(-x-y)/2}),
\end{displaymath}\nopagebreak
\begin{equation}
{\rm X}=\frac{1}{3}(e^{(y+z)/2}+e^{(y-z)/2}+e^{(-y-z)/2}), \label{mark1}
\end{equation}\nopagebreak
\begin{displaymath}
{\rm Y}=\frac{1}{3}(e^{(z+x)/2}+e^{(z-x)/2}+e^{(-z-x)/2}). 
\end{displaymath}

Using these expressions we can verify the equality
\begin{equation}
{\rm X}^2 + {\rm Y}^2 + {\rm Z}^2 - 3{\rm X}{\rm Y}{\rm Z} =
-\frac{1}{9}(e^{(x+y+z)/2} - e^{(-x-y-z)/2})^2
\end{equation}

The symmetry of the graph obviously cyclically permutes the coordinates and
therefore the numbers ${\rm Z},{\rm X},{\rm Y}$. A flip of an edge acts by
the rule (\ref{Hflipt}) and it results in the mapping

\begin{equation}
({\rm Z},{\rm X},{\rm Y}) \mapsto
({\rm Y},3{\rm Y}{\rm Z}-{\rm X},{\rm Z}).\label{flipmark}
\end{equation}

If all three coordinates $x,y,z$ are zeroes, the corresponding complex
surface is just the equiharmonic punctured torus.

The properties 1,3,5,6 immediately follows from this picture. One can easily
check, that $M(n)$ for $n \in {\Bbb N}$ are just the Fibonacci numbers what
gives the property 2. The property 7 is an immediate consequence of the
convexity property of the lamination length function. The property 4 was
proven by Markov himself.

The property 8 stands a little apart from the others since it is related to
the spaces ${\cal T}^d(S)$ and ${\cal L}^h(S)$ rather than ${\cal L}^d(S)$
and ${\cal T}^h(S)$, respectively. Consider a graph coordinate system $u,v,w$
on ${\cal T}^d(S)$; ${\rm U} = e^u,\ {\rm V}=e^v,\ {\rm W}=e^w$ and ${\rm A}$
is the area inside the horocycle. It easily follows from (\ref{area}) that
\begin{equation}
 ({\rm U}^2 + {\rm V}^2 + {\rm W}^2) = {\rm U}{\rm V}{\rm W}{\rm A}  \label{markdual}
\end{equation}
The cyclic symmetry of the graph acts by cyclic permutation of ${\rm U},{\rm V},{\rm W}$. A flip of an edge acts by
\begin{equation}
({\rm U},{\rm V},{\rm W}) \mapsto ({\rm W},
\frac{{\rm U}^2+{\rm W}^2}{\rm Z},{\rm U}).  \label{flipdecor}
\end{equation}

On the other hand this transformation law can be rewritten taking into account the equation (\ref{markdual}):
\begin{equation}
({\rm U},{\rm V},{\rm W}) \mapsto (({\rm W},{\rm U}{\rm W}{\rm A} - {\rm V},{\rm U})
\end{equation}

This rule coincides with (\ref{flipmark}) for $A=3$.

Now consider the decorated surface with ${\rm U}={\rm V}={\rm W}=1$. This is
the surface with the area inside the horocycle $A=3$. Applying modular
transformations we get obviously the Markov triples, what proves the
property 8.

There exists a canonical decomposition (called {\em main tesselation}) of
the upper half plane $H$ into ideal triangles with vertices in all rational
points of its ideal boundary. The dual graph to this tesselation is the
universal three-valent tree. The faces of this tree are therefore in
one-to-one correspondence with rational numbers. On the pictures below we
have drawn a fragment of this tree with corresponding Markov numbers written
on the faces.

\setlength{\unitlength}{0.25mm}
\begin{picture}(500,320)(100,500)
\put(400,780){\line( 0, 1){ 20}}
\put(400,800){\line( 1, 1){ 20}}
\put(420,820){\line( 1, 0){ 20}}
\put(400,800){\line(-1, 1){ 20}}
\put(420,820){\line(-1, 1){ 20}}
\put(240,700){\line( 4,-3){ 80}}
\put(320,640){\line( 1,-1){ 40}}
\put(360,600){\line( 1,-2){ 20}}
\put(160,640){\line( 1,-1){ 40}}
\put(200,600){\line( 1,-2){ 20}}
\put(320,640){\line(-1,-1){ 40}}
\put(280,600){\line(-1,-2){ 20}}
\put(280,600){\line( 1,-2){ 20}}
\put(360,600){\line(-1,-2){ 20}}
\put(200,600){\line(-1,-2){ 20}}
\put(120,600){\line( 1,-2){ 20}}
\put(400,780){\line(-2,-1){160}}
\put(240,700){\line(-4,-3){ 80}}
\put(160,640){\line(-1,-1){ 40}}
\put(120,600){\line(-1,-2){ 20}}
\put(400,780){\line( 2,-1){160}}
\put(560,700){\line( 4,-3){ 80}}
\put(640,640){\line( 1,-1){ 40}}
\put(680,600){\line( 1,-2){ 20}}
\put(680,600){\line(-1,-2){ 20}}
\put(640,640){\line(-1,-1){ 40}}
\put(600,600){\line(-1,-2){ 20}}
\put(600,600){\line( 1,-2){ 20}}
\put(560,700){\line(-4,-3){ 80}}
\put(480,640){\line(-1,-1){ 40}}
\put(440,600){\line(-1,-2){ 20}}
\put(440,600){\line( 1,-2){ 20}}
\put(480,640){\line( 1,-1){ 40}}
\put(520,600){\line( 1,-2){ 20}}
\put(520,600){\line(-1,-2){ 20}}
\put(100,560){\line(-2,-5){ 15.862}}
\put(100,560){\line( 2,-5){ 15.862}}
\put(140,560){\line(-2,-5){ 15.862}}
\put(140,560){\line( 2,-5){ 15.862}}
\put(180,560){\line(-2,-5){ 15.862}}
\put(180,560){\line( 2,-5){ 15.862}}
\put(220,560){\line(-2,-5){ 15.862}}
\put(220,560){\line( 2,-5){ 15.862}}
\put(260,560){\line(-2,-5){ 15.862}}
\put(260,560){\line( 2,-5){ 15.862}}
\put(300,560){\line(-2,-5){ 15.862}}
\put(300,560){\line( 2,-5){ 15.862}}
\put(340,560){\line(-2,-5){ 15.862}}
\put(340,560){\line( 2,-5){ 15.862}}
\put(380,560){\line(-2,-5){ 15.862}}
\put(380,560){\line( 2,-5){ 15.862}}
\put(420,560){\line(-2,-5){ 15.862}}
\put(420,560){\line( 2,-5){ 15.862}}
\put(460,560){\line(-2,-5){ 15.862}}
\put(460,560){\line( 2,-5){ 15.862}}
\put(500,560){\line(-2,-5){ 15.862}}
\put(500,560){\line( 2,-5){ 15.862}}
\put(540,560){\line(-2,-5){ 15.862}}
\put(540,560){\line( 2,-5){ 15.862}}
\put(580,560){\line(-2,-5){ 15.862}}
\put(580,560){\line( 2,-5){ 15.862}}
\put(620,560){\line(-2,-5){ 15.862}}
\put(620,560){\line( 2,-5){ 15.862}}
\put(660,560){\line(-2,-5){ 15.862}}
\put(660,560){\line( 2,-5){ 15.862}}
\put(700,560){\line(-2,-5){ 15.862}}
\put(700,560){\line( 2,-5){ 15.862}}
\small
\put(400,820){\makebox(0,0)[cb]{1}}
\put(400,760){\makebox(0,0)[cb]{5}}
\put(240,640){\makebox(0,0)[cb]{29}}
\put(560,640){\makebox(0,0)[cb]{13}}
\put(160,600){\makebox(0,0)[cb]{169}}
\put(320,600){\makebox(0,0)[cb]{433}}
\put(480,600){\makebox(0,0)[cb]{194}}
\put(640,600){\makebox(0,0)[cb]{34}}
\put(120,560){\makebox(0,0)[cb]{985}}
\put(200,560){\makebox(0,0)[cb]{14701}}
\put(280,560){\makebox(0,0)[cb]{37666}}
\put(360,560){\makebox(0,0)[cb]{6466}}
\put(440,560){\makebox(0,0)[cb]{2897}}
\put(520,560){\makebox(0,0)[cb]{7561}}
\put(600,560){\makebox(0,0)[cb]{1325}}
\put(680,560){\makebox(0,0)[cb]{89}}
\put(140,760){\makebox(0,0)[cb]{2}}
\put(660,760){\makebox(0,0)[cb]{1}}
\tiny
\put(100,500){\makebox(0,0)[cb]{5741}}
\put(140,500){\makebox(0,0)[cb]{499393}}
\put(180,500){\makebox(0,0)[cb]{7453378}}
\put(220,500){\makebox(0,0)[cb]{1278818}}
\put(260,500){\makebox(0,0)[cb]{3276569}}
\put(300,500){\makebox(0,0)[cb]{48928105}}
\put(345,500){\makebox(0,0)[cb]{8399329}}
\put(380,500){\makebox(0,0)[cb]{96557}}
\put(420,500){\makebox(0,0)[cb]{43261}}
\put(460,500){\makebox(0,0)[cb]{1686049}}
\put(500,500){\makebox(0,0)[cb]{4400489}}
\put(540,500){\makebox(0,0)[cb]{294685}}
\put(580,500){\makebox(0,0)[cb]{51641}}
\put(620,500){\makebox(0,0)[cb]{135137}}
\put(660,500){\makebox(0,0)[cb]{9077}}
\put(700,500){\makebox(0,0)[cb]{233}}
\end{picture}

\vspace{5mm}
\refstepcounter{equation}
\centerline{Fig. \arabic{equation}}

As a concluding remark of this section note that, as it was observed by
A.Bondal, Markov triples are dimensions of elements of distinguished sets of
sheaves on ${\Bbb C}P^2$. The relations between these two ways of obtaining
Markov numbers are completely unclear and very exciting.

\subsection{Duality between Teichm\" uller spaces.}

Here we are going to make some handwaving about what exactly we mean saying
that the Teichm\" uller spaces ${\cal T}^d$ and ${\cal T}^H$ are dual.  First
of all the meaning of duality between ${\cal L}^d$ and ${\cal L}^H$ can be
made precise. Indeed define the integral transform $S^{dh}:L^2({\cal L}^d)
\rightarrow L^2({\cal L}^h)$ given by

\begin{equation}
S^{dh}\psi(f_1) \mapsto \int_{{\cal L}^d} e^{i l_{\cal
LL}(f_1,f_2)}\psi(f_1){\rm vol}^d(f_1) \end{equation} where ${\rm vol}(f_1)$
is the canonical volume form on ${\cal L}^d$ given in graph coordinates
$\{u_\alpha\}$ by
\begin{equation} {\rm vol}^d = |\bigwedge_{\alpha \in
E(\Gamma)}du_\alpha|.
\end{equation}

Analogously one can define the conjugated integral
transform $S^{hd}:L^2({\cal L}^h) \rightarrow L^2({\cal L}^d)$ by
\begin{equation}
S^{hd}\psi(f_2) \mapsto \int_{{\cal L}^h} e^{i l_{\cal
LL}(f_1,f_2)}\psi(f_2){\rm vol}(f_2) \end{equation} where ${\rm vol}(f_1)$ is
the canonical volume form on ${\cal L}^h$. Define it first on ${\cal L}^H$ by
the analogous formula in graph coordinates $\{z_\alpha\}$:  \begin{equation}
{\rm vol}^H = |\bigwedge_{\alpha \in E(\Gamma)}dz_\alpha|.
\end{equation}

This volume form is also obviously invariant w.r.t. graph changes but also
w.r.t. changes of hole orientations thus defining the form ${\rm vol}^h$ on
${\cal L}^h$.

One can easily prove using the formula (\ref{LL}) and the Riemann
localisation theorem that the operators $S^{dh}$ and $S^{hd}$ are (up to a
scalar factor) mutually inverse isometries.

Of course, one can try to prove some analogous theorems for the integral
transforms related to the functions $l_{\cal LT}$ and $l_{\cal TL}$. However,
the most intriguing would be to formulate some analogous statements for the
Teichm\" uller spaces. There exist a construction of a function $l_{\cal TT}:
{\cal T}^d(S) \times {\cal T}^H(S) \rightarrow {\Bbb R}$ and asymptotically
compatible with the geodesic length functions $l_{\cal LT}$ and $l_{\cal TL}$
for closed surfaces due to F.Bonahon \cite{Bonahon}. It  is generalisable
for surfaces with holes, but unfortunately the Bonahon's construction is very
unexplicit and it is even very hard to prove that it gives a smooth function.
Of course it is hardly imaginable how one could check properties of the
corresponding integral transform without simplifying his approach. However,
as we have already mentioned in the introduction, it seems to us to be the
most important question in the whole subject to answer and we are going to do
it at least partially in the forthcoming preprint.

\subsection{Universal setting and Virasoro orbits.}

In this section we are going to use graph language to give some precise
sense to the statement that the simplest Virasoro coajoint orbit is a
universal Teichm\" uller space. The main idea of this notion belongs to
Bers (and was explained for me by R.C.Penner). Here we are just going to show
some simple statements that applying the graph language to the Bers
construction we get coordinates on the Virasoro coajoint orbit w.r.t. which
the canonical symplectic structure is constant.

 Call a {\em tesselation} a decomposition of a connected
part of the upper half plane $H$ in ideal triangles in such a way that each
ideal triangle has common edge with exactly three ideal triangles. The points
on the ideal boundary of $H$ being vertices of the triangles are called
vertices of the tesselation. The graph dual to the graph of edges of these
triangles is the universal three-valent fat tree. Call a tesselation
{\em full} if it covers $H$ completely. Denote by ${\it Tess}^h$ the space
of tesselations and by ${\it Tess}$ the space of full tesselations.

Particular examples of tesselations are given by the universal coverings of
Riemann surfaces with geodesic boundary cut into ideal triangles like in the
reconstruction procedure from section 5.1. Surfaces with all boundary
components being punctures give full tesselations.

 For any tesselation assign a positive real
number to each edge in a way analogous to what we have done for
graphs: take four vertices of the quadrilateral consisting of two triangles
separated by the edge, take a M\" obius function on $H$ taking values $0$ and
$\infty$ at the ends of the edge and $-1$ at a third vertex. The value $s$
of the M\" obius function at the fourth point is related to the number $z$ we
assign to the edge as $s=e^{z}$. Passing to the dual graph ${\tt T}$ one can
construct a tree with real numbers on edges out of any tesselation.

 Let ${\it Diff}$ and ${\it Homeo}$ be the groups of
smooth diffeomorphisms, and all
homeomorphisms of ${\Bbb R}P^1$, respectively. These groups act on the space of
tesselations ${\it Tess}^h$. Indeed, put into correspondence to each ideal triangle another
one vertices of which are images of the former one under the action of the
homeomorphism of the ideal boundary of $H$. This action is free and obviously
preserves the space ${\it Tess}$.

 One can check the following property of the space of tesselations:

\vspace{3mm}{\bf 1}. There exists one distinguished tesselation called {\em main
tesselation}. It is characterised by the properties that it contains the
ideal triangle with vertices at the points $0,1$ and $\infty$ and all its
graph coordinates are zeroes.

\vspace{3mm}{\bf 2}. The action of ${\it Homeo}$ on ${\it Tess}$ is transitive, therefore
we can identify the spaces ${\it Homeo}$ and ${\it Tess}$ by identifying the
identity with the main tesselation.

\vspace{3mm}{\bf 3}. Let ${\cal T_\infty} = {\it Tess}/PGL(2,{\Bbb R})$ be the space of
full tesselations modulo action of the standard $PGL(2,{\Bbb R})$ subgroup of
${\it Homeo}$. The association of numbers on edges on the tree to any
tesselation gives an embedding of ${\cal T}_\infty$ into ${\Bbb R}^\infty$.
The image consists of tesselations satisfying the following condition:
For any $\gamma \in F({\tt T})$ the sequences $0,e^{z_0}, e^{z_0 +
z_1} + e^{z_0}, \ldots$ and $-1, -e^{-z_{-1}}, -e^{-z_{-1}-z_{-2}} -
-e^{-z_{-1}}, \ldots$ diverge.  Here $\{z_i\}$ is a sequence of numbers
corresponding to edges belonging to $\gamma$ in their natural order.

\vspace{3mm}{\bf 4}. The 2-form $\omega_{WP}$ and the bivector $P_{WP}$ on ${\cal
T}_\infty$ given by
\begin{equation}
\omega_{WP} = \sum_{\aalpha \in EE({\tt T})}
dz_\aalpha \wedge dz_{\aalpha(1)}
\end{equation}
\begin{equation}
P_{WP} = \sum_{\aalpha \in EE(\Gamma)}\frac{\partial}{\partial z_\aalpha}
\wedge \frac{\partial}{\partial z_{\aalpha(1)}}
\end{equation}
are invariant
w.r.t. the action of ${\it Diff}$. The form $\omega_{WP}$ is closed. The
bivector $P_{WP}$ defines a Poisson structure on ${\cal T}_\infty$.

\vspace{3mm}{\bf 5}. Let ${\cal O} \subset {\cal T}_\infty$ be the image in ${\cal
T}_\infty$ of the group ${\it Diff}$. Then there exist a
${\it Diff}$-equivariant momentum mapping $\mu :{\cal O} \rightarrow
{\it vir}^*$, where ${\it vir}$ is the Virasoro algebra. The image of $\mu$
is the Virasoro coajoint orbit with the stabiliser $PGL(2,{\Bbb R})$.

\vspace{3mm}{\bf 6}. There exist canonical embeddings ${\cal T}(S) \rightarrow {\cal T}_\infty$ for any Riemann surface $S$. The image of a point of ${\cal T}(S)$ is given by the universal covering of a decomposition of the surface into
ideal triangles.

\section{Acknowledgments.} I am very grateful to the
  Max-Plank-Institut f\" ur Mathematik in Bonn and the Institut de Recherche
Math\' ematique Avanc\' ee in Strasbourg where most of the work has been
done. I am also very indebted to A.A.Rosly, R.C.Penner, Yu.Neretin,
A.Papadopulos, S.Orevkov and A.Zorich for valuable discussions, L.Shenderova
for drawing pictures, A.A.Vasiliev and H.Gangl for very attentive reading of the
manuscript and for A.Levin who explained me some applications of the
lamination theory.

The work was support in part by Award No. RM2-150 of the U.S. Civilian
Research and Development Foundation (CRDF) for the Independent States of
the Former Soviet Union.

The first part of this work was supported by the INTAS grant 96-518.

\end{document}